\RequirePackage{fix-cm} 
\documentclass[a4paper,twoside,12pt,reqno]{amsart}

\usepackage{fixltx2e}     

\usepackage[english]{babel}
\usepackage[latin1]{inputenc}

\usepackage{indentfirst,verbatim}

\usepackage{amsmath,amsfonts,amssymb,amsgen,amsbsy,eucal,mathrsfs,dsfont}
\usepackage{stmaryrd}

\usepackage[square,numbers]{natbib}
\usepackage{a4wide,verbatim}

\usepackage{epsfig,rotating,color}
\usepackage{multicol}
\usepackage{tikz,pgfplots}
\usepackage{graphicx}
\usepackage{mathtools}
\usepackage{textcomp}

\usepackage{caption}

\usepackage{hyperref}

\usepackage[square,numbers]{natbib}
\usepackage{systeme}

\hypersetup{
colorlinks=true, 
breaklinks=true, 
urlcolor= blue, 
linkcolor= blue, 
bookmarksopen=true, 
pdftitle={}, 
pdfauthor={}, 
pdfsubject={} 
}

\newcommand{\R}{\mathds{R}}

\newcommand{\ud}{\mathrm{d}}

\newcommand{\half}{{\textstyle{1\over2}}}

\usepackage{amsthm}
\newtheorem{thm}{Theorem}[section]
\newtheorem{definition}[thm]{Definition}
\newtheorem{lem}[thm]{Lemma}

\theoremstyle{remark}
\newtheorem{rem}[thm]{Remark}

\usepackage{enumitem}
\newlist{steps}{enumerate}{1}
\setlist[steps, 1]{label = Step \arabic*:}

\newcommand{\eqdef}{\stackrel{\text{\tiny{def}}}{=}}

\title[Regularized scalar conservation laws]{\bf On a Hamiltonian regularization of scalar conservation laws}

\author[GUELMAME]{Billel Guelmame}

\newcommand{\nfont}{\fontshape{n}\selectfont}

\address{({\nfont\textbf{Billel Guelmame}})  LJAD,  Inria, UMPA, CNRS,  Universit\'e C\^ote d'Azur, ENS Lyon, France.} 
\email{billel.guelmame@univ-cotedazur.fr}

\newcommand\smallO{
  \mathchoice
    {{\scriptstyle\mathcal{O}}}
    {{\scriptstyle\mathcal{O}}}
    {{\scriptscriptstyle\mathcal{O}}}
    {\scalebox{.7}{$\scriptscriptstyle\mathcal{O}$}}
  }

\let\oldtocsection=\tocsection
 
\let\oldtocsubsection=\tocsubsection

\renewcommand{\tocsection}[2]{\hspace{0em}\oldtocsection {#1}{#2}}
\renewcommand{\tocsubsection}[2]{\hspace{2em}\oldtocsubsection{#1}{#2}}

\numberwithin{equation}{section}

\usepackage{tabularx}
\usepackage[pagewise]{lineno}  

\begin{document}

\maketitle

\begin{abstract}
In this paper, we propose a Hamiltonian regularization of scalar conservation laws, which is parametrized by $\ell>0$ and conserves an $H^1$ energy. We prove the existence of global weak solutions for this regularization. Furthermore, we demonstrate that as $\ell$ approaches zero, the unique entropy solution of the original scalar conservation law is recovered, providing justification for the regularization.

This regularization belongs to a family of non-diffusive, non-dispersive regularizations that were initially developed for the shallow-water system and  extended later to the Euler system.
This paper represents a validation of this family of regularizations in the scalar case.
\end{abstract}

\medskip

 {\bf AMS Classification :} 35L65; 35B65; 35L67; 35Q35.

\medskip

{\bf Key words :} Scalar conservation laws; Hamiltonian; regularization; 
conservative and dissipative solutions; Oleinik inequality.

\tableofcontents

\section{Introduction}

Hyperbolic conservation laws, such as the inviscid Burgers equation and the barotropic Euler system, are known to develop discontinuous shocks even when the initial data is a smooth ${C}^\infty$ function. This poses a challenge both in numerical simulations and in theoretical studies. In order to avoid these discontinuous shocks, diffusion and/or dispersion terms can be added into the equations.
In \cite{ClamondDutykh2018a}, Clamond and Dutykh derived a non-diffusive, non-dispersive regularized Saint-Venant (rSV) system, which is Galilean invariant and conserves an $H^1$-like energy for smooth solutions.
The weakly singular shock profiles of the rSV system have been studied in \cite{PuEtAl2018}, while the local well-posedness and the blow-up scenarios for the rSV system have been studied in \cite{LiuEtAl2019}.
A regularized barotropic Euler (rE) system was proposed and studied in \cite{guelmame2020Euler} as a generalization of the rSV system.
Both the rSV and rE systems are locally well-posed in $H^s$ with $s \geqslant 2$, however, their solutions may develop singularities in finite time \cite{LiuEtAl2019, SGNbu2021}. 
The study of these systems remains challenging, with both the existence of global weak solutions and the understanding of the limiting case being outstanding problems.
However, due to the similarities between the rSV, rE, and the Serre--Green--Naghdi (SGN) systems, it may be possible to obtain global weak solutions for small-data of the rSV and rE systems following a recent proof for the SGN with surface tension \cite{guelmameSGN2}.

Inspired by the rSV system, the regularized Burgers (rB) equation 
\begin{equation}\label{rB0}
u_t\ +\ u\, u_x\ =\ \ell^2\, \left[u_{txx}\ +\ 2\, u_x\, u_{xx}\, +\ u\, u_{xxx} \right]
\end{equation} 
have been proposed in \cite{guelmamerB},
where $\ell$ is a positive parameter.
Being a scalar equation, the rB equation is more tractable than the rSV system.
In \cite{guelmamerB}, a study of weakly singular shocks and cusped traveling-wave weak solutions of \eqref{rB0} is established. 
Also, inspired by \cite{BressanConstantin2007a,BressanConstantin2007b}, a proof the existence of two types of global weak solutions of \eqref{rB0}, conserving or dissipating the energy is presented.
The dissipative solution of \eqref{rB0} satisfies the one-sided Oleinik inequality
\begin{equation}\label{Ol0}
 u_x(t,x)\ \leqslant\ 2/t, \qquad \forall\, (t,x)\ \in\ (0,\infty) \times \R.
\end{equation}
The limiting cases of $\ell \to 0$ and $\ell \to \infty$ have also been studied in \cite{guelmamerB}. However, the equations satisfied by the limits were not well established.
The rB equation \eqref{rB0} must be compared to the well-known dispersionless Camassa--Holm (CH) equation \cite{CamassaHolm1993}
\begin{equation}\label{CH}
u_t\ +\ 3\, u\, u_x\ =\ \ell^2\, \left[u_{txx}\ +\ 2\, u_x\, u_{xx}\, +\ u\, u_{xxx} \right].
\end{equation} 
Both \eqref{rB0} and \eqref{CH} conserve an $H^1$ energy (not uniformly on $\ell$) for smooth solutions.
A key difference between the two equations is that \eqref{rB0} is Galilean invariant while \eqref{CH} is not. This Galilean invariance is significant not only from a physical point of view, but also mathematically. Due to the Galilean invariance of the rB equation, the constant on the right-hand side of the Oleinik inequality \eqref{Ol0} is independent of $\ell$. On the other hand, the dissipative solutions of the CH equation satisfy a similar inequality as \eqref{Ol0}, but with a constant that depends on $\ell$. As a result, the compactness arguments presented in \cite{guelmamerB} cannot be used for the CH equation.
However, the limiting case of the viscous CH equation
have been studied in \cite{CD2016,CK1,Hwang2007} under the condition ``$\ell$ is small enough compared to the
viscosity parameter''. The authors proved that as the viscosity parameter goes to zero, the unique entropy solution of the scalar conservation law $u_t+ (3u/2)_x=0$ is recovered.

This paper is a continuation of the previous one \cite{guelmamerB}. Our goal is to generalize the rB equation \eqref{rB0} to regularize scalar conservation laws, to prove the existence of global weak solution of the regularized equation, and to study the limiting cases $\ell \to 0$ and $\ell \to \infty$.
We consider the equation 
\begin{equation}\label{rSCL}
u_t\ +\ f(u)_x\ =\ \ell^2\, \left[ u_{xxt}\ +f'(u)\, u_{xxx}\ +\ 2\, f''(u)\, u_x\, u_{xx}\  +\ \half\,  f'''(u)\, u_x^3 \right],
\end{equation}
where $f$ is a uniformly convex ($f''(u) \geqslant c > 0$) flux,
the rB equation \eqref{rB0} is recovered taking $f(u) = u^2 /2$.
The equation \eqref{rSCL} has several interesting properties, such as conservation of an $H^1$ energy for smooth solutions, and both Hamiltonian and Lagrangian structures. 
Holden and Raynaud \cite{HR20072} conducted a study on a generalized version of both the Camassa--Holm equation and the hyperelastic-rod wave equation. This study includes the equation \eqref{rSCL}, which we derived here from a distinct motivation to regularize scalar conservation laws.
The existence of global weak solutions of the Camassa--Holm equation and its generalizations in the space $H^1$ has been widely studied before.
There are two types of solutions, conserving and dissipating the energy.
The proof of the existence of conservative solutions uses equivalent systems of ODEs written in the Lagrangian coordinates \cite{BressanConstantin2007a,HR20072}.
Conservative solutions fail to satisfy the one-sided Oleinik inequality \cite{guelmamerB}, which is a crucial property of entropy solutions of scalar conservation laws. Hence, to regularize scalar conservation laws, we need to consider dissipative solutions of \eqref{rSCL}.
Dissipative solutions can be obtained through various methods, such as equivalent systems in the Lagrangian coordinates  \cite{BressanConstantin2007b}, vanishing viscosity \cite{CHK,XinZhang2000}, and the convergence of finite difference
schemes \cite{CKR4,CKR5}.
In this paper, we demonstrate the existence of global weak solutions of \eqref{rSCL} with a different method. 
Our approach involves an approximated equation through a cut-off in the Riccati equation, similar to the methods employed in \cite{wave3}.
Our approximated equation is globally well-posed, we obtain then some uniform estimates that allow us to use classical compactness arguments with Young measures \cite{Focusing}. Taking the limit in the approximated equation leads to the global dissipative solution of \eqref{rSCL}.

The formal limit  $\ell \to 0$ and $\ell \to \infty$ in \eqref{rSCL} lead to the equations 
\begin{subequations}
\begin{align}\label{SCL}
u_t\ +\ f(u)_x\ &=\ 0 \qquad & \mathrm{as}\ \ell\ \to\ 0,\\ \label{gHS}
\left[ u_t\ +\ f(u)_x \right]_x\, &=\ \half\, u_x^2\, f''(u) & \mathrm{as}\ \ell\ \to\ \infty.
\end{align}
\end{subequations}
Equation \eqref{SCL} is known as the scalar conservation law, while \eqref{gHS} is a generalized Hunter--Saxton equation \cite{HunterSaxton1991}.
The classical Hunter--Saxton (HS) equation is recovered taking $f(u) = u^2/2$.
The HS equation and its generalization \eqref{gHS} admit conservative and dissipative global weak solutions \cite{BressanConstantin05,BZ2007,HunterZheng1995}.
In a previous study \cite{guelmamerB}, the limiting cases $\ell \to 0$ and $\ell \to \infty$ of \eqref{rB0} were investigated. However, the two limits satisfy equations that involve Radon measures which were not identified.
This present work proves that, when $\ell \to 0$, the dissipative solution of \eqref{rSCL} converges to the unique entropy solution of \eqref{SCL}. Furthermore, as $\ell \to \infty$, the dissipative solution of \eqref{rSCL} converges to a dissipative solution of the generalized Hunter--Saxton equation \eqref{gHS}.
In the case where $\ell \to 0$, we present improved estimates compared to those found in \cite{guelmamerB}. These improved estimates enable us to identify the Radon measure that was left unidentified in \cite{guelmamerB}.
In order to identify the Radon measure in the case of $\ell \to \infty$, we employ some Young measures techniques. This enables us to take the limit of a quadratic term while only having a weak limit.

This paper is organized as follows. In Section \ref{sec:properties}, we present the variational formulations and other properties of the equation \eqref{rSCL}. We state the main results of the paper in Section \ref{sec:mr}. 
Section \ref{sec:ae} presents the approximated equation of \eqref{rSCL} and provides some uniform estimates.
Section \ref{sec:precompactness} establishes the existence of a global dissipative solution.
The limiting case $\ell \to 0$ is studied in Section \ref{sec:zero}. Finally, Section \ref{sec:infty} focuses on the case $\ell \to \infty$.

\section{Variational formulations}\label{sec:properties}
This section is devoted to present some properties of the regularized equation \eqref{rSCL}, including its Hamiltonian and Lagrangian structures.

Applying the operator $\left(\,1 -\ell^2 \partial_x^{\,2}\,\right)^{\!-1}$ to \eqref{rSCL}, we obtain 
\begin{equation}\label{rSCLinv}
u_t\ +\ f(u)_x\ +\ \half\,\ell^2\left(\,1\,-\,\ell^2\,
\partial_x^{\,2}\,\right)^{\!-1}\left[\, f''(u)\ u_x^{\,2}\,\right]_x\ =\ 0.
\end{equation}
The equation \eqref{rSCL} can be obtained as the Euler--Lagrange equation of the Lagrangian density
\begin{equation*}
\mathcal{L}_\ell\ \eqdef\ \half\,\phi_x\,\phi_t\ +\ F(\phi_x) +\ 
\half\, \ell^2\ \left[f'(\phi_x)\ \phi_{xx}^2\ - \phi_{xxx}\phi_t \right], \qquad \phi_x\ =\ u,
\end{equation*}
where $F'(u)=f(u)$.
A Hamiltonian structure also exists for the equation \eqref{rSCLinv}, that can be obtained with the Hamiltonian operator and functional 
\begin{align}\label{Hamiltonian}
\mathscr{D}\ \eqdef\ \left(1\ -\ \ell^2\, \partial_x^2 \right)^{-1}\,\partial_x, \qquad \quad
\mathfrak{H}\ \eqdef\ \int\left[\, F(u)\,+\,
\half\,\ell^2\, f'(u)\ u_x^2\, \right]\,\ud\/x,
\end{align}
so the equation of motion is given by 
\begin{align*}
u_t\ =\ -\,\mathscr{D}\ \delta_u\,\mathfrak{H},
\end{align*}
where the operator $\mathscr{D}$ is a Hamiltonian operator. 
Defining
\begin{equation}\label{Pdef}
P\ \eqdef\ \half \left(\,1 -\ell^2 \partial_x^{\,2}\,\right)^{\!-1} \left\{ f''(u)\, u_x^{2}\right\}\, =\ \half\, \mathfrak{G}\ast \left\{ f''(u)\, u_x^{2}\right\},
\end{equation}
where
\begin{equation*}
\mathfrak{G}\ \eqdef\ (2\ell)^{-1}\,\exp(-|\cdot|/\ell).
\end{equation*}
Smooth solutions of \eqref{rSCLinv} satisfy the energy conservation
\begin{equation}
\left[\,\half\,u^2\ +\ \half\,\ell^2\,u_x^{\,2}\,\right]_t\ 
+\ \left[\, K(u)\ +\ \half\, \ell^2\, f'(u)\, u_x^2\ +\ \ell^2\, u\, P\, \right]_x\ =\ 0,
\end{equation}
where $K'(u)=uf'(u)$. Another conservation equation that corresponds to the Hamiltonian \eqref{Hamiltonian} can be obtained
\begin{equation*}
\left[F(u)\, +\, \half\,\ell^2\, f'(u)\ u_x^2 \right]_t\, +\, \left[\half\, f(u)^2\, +\, \ell^2\, f(u)\, P\, +\, \half\, \ell^2\, f'(u)^2\, u_x^2\, +\, \half\, \ell^4\, P^2\, -\, \half\, \ell^6\, P_x^2 \right]_x\, =\ 0.
\end{equation*}
In next section, we present the main results of this paper.

\section{Main results}\label{sec:mr}
We consider the Cauchy problem
\begin{equation}\label{rB}
u^\ell_t\ +\,\left[\,f(u^\ell)\,+\,\half\,\ell^2\,\mathfrak{G}\ast \left\{ f''(u^\ell) \left(u^\ell_x\right)^{2}\right\}
\right]_x\ =\ 0,
\end{equation}
with $u^\ell(0,x) = u_0(x)$. 
Using that $P - \ell^2 P_{xx} =  f''(u^\ell) \left(u^\ell_x\right)^2/2$ and
differentiating \eqref{rB} w.r.t $x$ one obtains
\begin{equation}\label{barqeqell}
u^\ell_{xt}\, +\, \left[f'(u^\ell)\, u_x^\ell\right]_{x}\  =\ -\ P\, +\ \half\, f''(u^\ell) \left(u_x^\ell\right)^2.
\end{equation}
We start this section by defining dissipative weak solutions of \eqref{rB}.
\begin{definition}\label{WSDef} We say that $u^\ell \in L^\infty(\mathds{R}^+, H^1) \cap \mathrm{Lip}(\mathds{R}^+, L^2)$ is a weak dissipative solution of \eqref{rB} if it satisfies the initial condition $u^\ell(0,\cdot)=u_0$ with \eqref{rB} in the $L^2$ sense and dissipates the energy
\begin{gather}\label{Ene_local}
\left[\,\half \left(u^{\ell}\right)^2\ +\ \half\,\ell^2 \left(u_x^{\ell}\right)^{\,2}\,\right]_t\ 
+\ \left[\, K \left(u^{\ell}\right)\ +\ \half\, \ell^2\, f'\! \left(u^{\ell}\right)\, \left(u_x^{\ell}\right)^2\ +\ \ell^2\, u^{\ell}\, P\, \right]_x\, \leqslant\ 0.
\end{gather}
Moreover, $u^\ell$ is right continuous in $H^1$. More precisely, for all $t_0 \geqslant 0$ we have
\begin{equation}\label{rcont}
\lim_{\substack{t\to t_0\\ t> t_0}} \left\| u^\ell(t,\cdot)\, -\, u^\ell(t_0,\cdot)\right\|_{H^1}\ =\ 0.
\end{equation}
\end{definition}

\begin{thm}\label{thm:existence}
Let $f \in C^4$ be a uniformly convex flux ($f''(u) \geqslant c >0$),  $u_0 \in H^1(\R)$ and $\ell>0$, then there exists a global weak dissipative solution $u^\ell \in L^\infty ([0,\infty ), H^1(\mathds{R})) \cap C([0,\infty ) \times \mathds{R})$ of \eqref{rB} in the sense of Definition \ref{WSDef} satisfying the following
\begin{itemize}
\item For any $T>0$, any bounded set $ [a,b] \subset \mathds{R}$ and $\alpha \in [0,1)$ there exists $C=C(\alpha, T, a, b, \ell)>0$ such that 
\begin{equation}\label{alpha+2_xi}
\int_0^T \int_a^b \left[\left|u^\ell_t\right|^{2+\alpha}\ +\ \left|u^\ell_x\right|^{2+\alpha} \right] \mathrm{d}x\, \mathrm{d}t\ \leqslant\ C . 
\end{equation}
\item The solution satisfies the one-sided Oleinik inequality
\begin{equation}\label{Ol:main:thm}
u_x^{\ell}(t,x)\ \leqslant\ \frac{1}{c\, t/2\, +\, 1/M} \qquad a.e.\ (t,x) \in (0,\infty) \times \R,
\end{equation}
where $M=\sup_x u_0'(x) \in (0,\infty]$. 
\end{itemize}
Moreover, if $f''(u) \leqslant C$, $u_0' \in L^1(\R)$ and $u_0'(x) \leqslant M < \infty$ then 
\begin{equation}\label{TV}
\left\| u^{\ell} \right\|_{L^\infty}\, \leqslant\, \left\| u_x^{\ell} \right\|_{L^1}\, \leqslant\ \left\| u_0' \right\|_{L^1} \left( c\, M\, t/2\, +\, 1 \right)^{2\, C/c}, \qquad \forall \ell \in (0,\infty).
\end{equation}
\end{thm}

\begin{rem}
The constant $C>0$ in \eqref{alpha+2_xi} depends also on $\ell$. 
In Lemma \ref{lem:alpha+22} below, we prove that if $\ell$ is far from $0$, one can chose a constant $C$ independent on $\ell \geqslant 1$.
\end{rem}

The aim of this paper is to prove that the equation \eqref{rB} is indeed a regularisation of scalar conservation laws. i.e., as $\ell \to 0$ the dissipative solution of \eqref{rB} giving by Theorem \ref{thm:existence} converges to the unique entropy solution of the scalar conservation law \eqref{SCL}.
\begin{thm} \label{thm:0}
Let $f \in C^4$ be a uniformly convex flux such that $C \geqslant f''(u) \geqslant c >0$. Let $u_0 \in H^1(\R)$ such that $u_0' \in L^1(\R)$ and  $u_0'(x) \leqslant M < \infty$, and let also $u^\ell$ be the dissipative solution of \eqref{rB} given by Theorem \ref{thm:existence}, then there exists a limit $u^0 \in L^\infty_{loc}([0,\infty),L^1_{loc}(\R))$ such that
\begin{itemize}
\item $  u^\ell \xrightarrow{\ell \to 0} u^0$ in $L^\infty_{loc}([0,\infty),L^p_{loc}(\R))$ for all $p \in [1,\infty)$.
\item $u^0$ is the unique entropy solution of the scalar conservation law \eqref{SCL}.
\end{itemize}
\end{thm}

As mentioned above, taking $\ell \to \infty$ formally in \eqref{rB} we obtain the generalized Hunter--Saxton equation \eqref{gHS}. 
We prove here that, up to a subsequence, the dissipative solution of  \eqref{rB} converges to a dissipative solution of \eqref{gHS} as $\ell \to \infty$.

\begin{definition}\label{WSDefHS} We say that $u \in L^\infty(\mathds{R}^+, \dot{H}^1(\R))$ is a weak dissipative solution of \eqref{gHS} if it satisfies the initial condition $u^\ell(0,\cdot)=u_0$ with \eqref{gHS} in the sense of distributions and dissipates the energy
\begin{gather}\label{Ene_localHS}
\left[ \left(u_x^{\infty}\right)^{\,2}\,\right]_t\ 
+\ \left[ f'\! \left(u^{\infty}\right)\, \left(u_x^{\infty}\right)^2 \right]_x\, \leqslant\ 0.
\end{gather}
Moreover, $u$ is right continuous in $\dot{H}^1$. More precisely, for all $t_0 \geqslant 0$ we have
\begin{equation}\label{rcont2}
\lim_{\substack{t\to t_0\\ t> t_0}} \left\| u(t,\cdot)\, -\, u(t_0,\cdot)\right\|_{\dot{H}^1}\ =\ 0.
\end{equation}
\end{definition}

\begin{thm} \label{thm:inf}
Let $f \in C^4$ be a uniformly convex flux ($f''(u) \geqslant c >0$),  $u_0 \in H^1(\R)$  and let $u^\ell$ be the dissipative solution of \eqref{rB} given by Theorem \ref{thm:existence}, then there exists a subsequence of $\left(u^\ell\right)_\ell$ that we denote also $\left(u^\ell\right)_\ell$ and a limit $u^\infty \in L^\infty_{loc}([0,\infty)\times \R) \cap L^\infty([0,\infty),\dot{H}^1(\R))$ such that
\begin{itemize}
\item $  u^\ell \xrightarrow{\ell \to \infty} u^\infty$ in $L^\infty_{loc}([0,\infty)\times \R) \cap \dot{H}^{1}_{loc}([0,\infty)\times \R)$.
\item $u^\infty$ is a dissipative solution of the generalized Hunter--Saxton equation \eqref{gHS}.
\item 
$u_x^{\infty}(t,x)\ \leqslant\ \frac{1}{c\, t/2\, +\, 1/M}$ $a.e.\ (t,x) \in (0,\infty) \times \R,$
where $M=\sup_x u_0'(x) \in (0,\infty]$.
\item $u^\infty_t,u^\infty_x \in L^{2+\alpha}_{loc}([0,\infty)\times \R)$, $\forall \alpha \in [0,1)$.
\end{itemize}
Moreover, if $f''(u) \leqslant C$, $u_0' \in L^1(\R)$ and $u_0'(x) \leqslant M < \infty$ then 
\begin{equation}\label{TVgHS}
\left\| u_x^{\infty} \right\|_{L^1}\, \leqslant\ \left\| u_0' \right\|_{L^1} \left( c\, M\, t/2\, +\, 1 \right)^{2\, C/c}.
\end{equation}
\end{thm}

\begin{rem} \hfill 
\begin{itemize}
\item If $f(u)=u^2/2$, the equation \eqref{gHS} is the classical Hunter--Saxton equation, and $u^\infty$ is the unique \cite{Dafermos} dissipative solution of \eqref{gHS}.
\item The proof presented in this paper of the limiting case $\ell \to 0$ (Theorem \ref{thm:0}) cannot be used for the Camassa--Holm equation \eqref{CH}.
\item The proof of Theorem \ref{thm:inf} (except \eqref{TVgHS}) works for the Camassa--Holm equation. In other words, the dissipative solutions of the Camassa--Holm equation \eqref{CH} converge to the dissipative solutions of the Hunter--Saxton equation (Eq. \ref{gHS} with $f(u)=u^2/2$) as $\ell \to \infty$.
\end{itemize}
\end{rem}

\section{The approximated equation}\label{sec:ae}

The local (in time) well-posedness of the equation \eqref{rB} can be obtained easily using Kato's theorem for quasi-linear equations \cite{Katoquasi}. However, for uniformly convex fluxes non-trivial solutions blow-up in a finite time. The blow-up occurs due the quadratic terms in \eqref{barqeqell}, this can be proved simply by following the characteristics in the equation \eqref{barqeqell} and using that $P\geqslant 0$. 
In order to prove the existence of global dissipative solutions of \eqref{rB}, we use a cut-off function instead of the quadratic terms of \eqref{barqeqell} that behaves like a linear function when $u_x^\ell$ is very large and negative (near to $-\infty$).
This cut-off is the key point to obtain global solutions.
As in \cite{wave3}, we define for any $\varepsilon>0$
\begin{equation}\label{chidef}
\chi_\varepsilon (q)\ \eqdef\, \left(q\ +\ \frac{1}{\varepsilon} \right)^2 \mathds{1}_{(-\infty,-\frac{1}{\varepsilon}]} (q)\ =\ 
\begin{cases}
\left(q\ +\ \frac{1}{\varepsilon} \right)^2, & q \leqslant -1/\varepsilon, \\
0, & q > -1/\varepsilon.
\end{cases}
\end{equation}
In order to obtain smooth solutions of the truncated equation, we choose $\chi_\varepsilon$ to be $C^1$ and piecewise $C^2$ instead of a rough cut-off.

Let $u_0 \in H^1$ and $j_\varepsilon$ be a Friedrichs mollifier, we define $u_0^\varepsilon \eqdef u_0 \ast j_\varepsilon$ and we consider the approximated Cauchy problem
\begin{equation}\label{rBep}
u^{\ell,\varepsilon}_t\ +\, \left[f(u^{\ell,\varepsilon})\,+\,\half\,\ell^2\,\mathfrak{G}\ast \left\{ f''(u^{\ell,\varepsilon}) \left\{ (u^{\ell,\varepsilon}_x)^{\,2} +\, \chi_\varepsilon(u_x^{\ell,\varepsilon}) \right\} \right\}
\right]_x\, =\ 0, \quad u^{\ell,\varepsilon}(0,\cdot)\, =\, u_0^\varepsilon.
\end{equation}
Defining 
\begin{equation}\label{Pepdef}
q^{\ell,\varepsilon}\ \eqdef\ u^{\ell,\varepsilon}_x, \qquad P^\varepsilon\ \eqdef\ \half\,\mathfrak{G}\ast \left\{ f''(u^{\ell,\varepsilon}) \left\{ (u^{\ell,\varepsilon}_x)^{\,2} +\, \chi_\varepsilon(u_x^{\ell,\varepsilon}) \right\} \right\}.
\end{equation}
Differentiating \eqref{rBep} with respect to $x$ we obtain 
\begin{equation}\label{rBep_x}
q^{\ell,\varepsilon}_{t}\, +\ f'(u^{\ell,\varepsilon})\, q^{\ell,\varepsilon}_{x}\ +\ \half\, f''(u^{\ell,\varepsilon})\, (q^{\ell,\varepsilon})^2\ +\ P^\varepsilon\ -\ \half\, f''(u^{\ell,\varepsilon})\, \chi_\varepsilon(q^{\ell,\varepsilon})\ =\ 0.
\end{equation}
Multiplying \eqref{rBep} by $u^{\ell,\varepsilon}$ and \eqref{rBep_x} by $\ell^2 q^{\ell,\varepsilon}$ we obtain 
the energy equation
\begin{gather} \nonumber
\left[\,\half \left(u^{\ell,\varepsilon}\right)^2\ +\ \half\,\ell^2 \left(q^{\ell,\varepsilon}\right)^{\,2}\,\right]_t\ 
+\ \left[\, K \left(u^{\ell,\varepsilon}\right)\ +\ \half\, \ell^2\, f'\! \left(u^{\ell,\varepsilon}\right)\, \left(q^{\ell,\varepsilon}\right)^2\ +\ \ell^2\, u^{\ell,\varepsilon}\, P^\varepsilon\, \right]_x\\ \label{Ene_ep}
 =\ \half\, \ell^2\, f''(u^{\ell,\varepsilon})\, q^{\ell,\varepsilon}\, \chi_\varepsilon(q^{\ell,\varepsilon})\ \leqslant\ 0.
\end{gather}
Our goal is to prove that the approximated equation \eqref{rBep} admits global smooth solutions, and, taking $\varepsilon \to 0$ we obtain global weak solutions of \eqref{rB}. We present now the existence of global solutions of \eqref{rBep}.
\begin{thm}\label{Thm:existenceapp}
Let $f \in C^4$, $\ell,\varepsilon>0$ and $u_0 \in H^1$, there exists a global smooth solution $u^{\ell,\varepsilon} \in C(\mathds{R}^+,H^3(\mathds{R})) \cap C^1(\mathds{R}^+,H^2(\mathds{R}))$ of \eqref{rBep} satisfying \eqref{Ene_ep}.
\end{thm}
The proof of Theorem \ref{Thm:existenceapp} is classical and is omitted here. 
The local well-posedness of \eqref{rBep} can be obtained using Kato's theorem for quasi-linear hyperbolic equations \cite{Katoquasi}. 
When $q^{\ell,\varepsilon} \leqslant -1/\varepsilon$, the quadratic term in the Riccati equation \eqref{rBep_x} becomes linear, this prevents the singularities from appearing in finite time and leads to the global well-posedness of \eqref{rBep}.
\qed 

Integrating \eqref{Ene_ep}, we obtain
\begin{gather}\nonumber
\int_\R \left[ (u^{\ell,\varepsilon})^2\, +\, \ell^2\, (u_x^{\ell,\varepsilon})^2 \right] \ud x\ -\ \ell^2 \int_0^t \int_\R f''(u^{\ell,\varepsilon})\, u_x^{\ell,\varepsilon}\, \chi_\varepsilon(u_x^{\ell,\varepsilon})\, \ud x\, \ud t\\ \label{eneequation}
=\ \int_\R \left[ (u_0^\varepsilon)^2\, +\, \ell^2\, \left(\partial_x u_0^\varepsilon\right)^2 \right] \ud x\
\leqslant\ \int_\R \left[ (u_0)^2\, +\, \ell^2\, \left( u_0'\right)^2 \right] \ud x.
\end{gather}
The energy equation \eqref{eneequation} implies that 
\begin{subequations}
\begin{gather}
\|u^{\ell,\varepsilon}\|_{L^2}\ \leqslant\, \left( \ell^2\, +\, 1 \right)^\frac{1}{2} \|u_0\|_{H^1},\\ \label{uxenergy}
\|u_x^{\ell,\varepsilon}\|_{L^2}\ \leqslant\, \left( \ell^{-2}\, +\, 1 \right)^\frac{1}{2} \|u_0\|_{H^1}.
\end{gather}
\end{subequations}
Then, the embedding $H^1 \hookrightarrow L^\infty$ implies that  
\begin{gather}\label{uLinf}
\|u^{\ell,\varepsilon}\|_{L^\infty}\ \leqslant\ C_\ell\, \|u_0\|_{H^1}.
\end{gather}
Using that $\chi_\varepsilon(q) \leqslant q^2$ and Young inequality we obtain for all $p \in [1,\infty]$ that
\begin{subequations}\label{PLp}
\begin{gather}
\|P^\varepsilon\|_{L^p}\ \leqslant\ C\, \|\mathfrak{G}\|_{L^p}\, \left\|u_x^{\ell,\varepsilon} \right\|_{L^2}^2\ \leqslant\ \tilde{C}\, \ell^\frac{1-p}{p} \left( \ell^{-2}\, +\, 1 \right),\\
\|P_x^\varepsilon\|_{L^p}\ \leqslant\ C\, \|\mathfrak{G}_x\|_{L^p}\, \left\|u_x^{\ell,\varepsilon} \right\|_{L^2}^2\ \leqslant\ \tilde{C}\, \ell^\frac{1-2p}{p} \left( \ell^{-2}\, +\, 1 \right).
\end{gather}
\end{subequations}
 We prove now a one-sided Oleinik inequality that is uniform on $\ell>0$ and $\varepsilon>0$. This inequality is very important. It is used in Lemma \ref{lemTV} below to obtain a uniform estimate of the total variation of the solution $u^{\ell,\varepsilon}$. Moreover, it is used in Sections \ref{sec:precompactness} and \ref{sec:infty} below to study the limiting cases $\varepsilon \to 0$ and $\ell \to \infty$.
It is also used in Section \ref{sec:zero} to show that the limit as $\ell \to 0$ is the unique entropy solution of the scalar conservation law \eqref{SCL}.
\begin{lem}\label{Lem:Oleinik}
Let $u_0 \in H^1$, $f''(u) \geqslant c>0$ and $u^{\ell,\varepsilon}$ be the solution of \eqref{rBep} given by Theorem \ref{Thm:existenceapp}, then 
\begin{equation}\label{Oleinikellep}
u_x^{\ell,\varepsilon}(t,x)\ \leqslant\ \frac{1}{c\, t/2\, +\, 1/M} \qquad a.e.\ (t,x) \in (0,\infty) \times \R,
\end{equation}
where $M=\sup_x u_0' \in [0,\infty]$.
\end{lem}
\proof For a fixed $x \in \R$, let $X(\cdot,x)$ be the solution of the ODE $X_t(t,x) = u^{\ell, \varepsilon}(t,X(t,x))$ with $X(0,x)=x$. Let $h(t) \eqdef u^{\ell,\varepsilon}_x(t,X(t,x))$, the equation \eqref{rBep_x} implies that 
\begin{equation}
h'(t)\ \leqslant\ - \half\, c\, h(t)^2\ +\ \half\, f''(u^{\ell, \varepsilon}(t,X(t,x)))\, \chi_\varepsilon(h(t)).
\end{equation}
Initially, $h(0)=\partial_x u_0^\varepsilon (x) \leqslant \sup_x u_0'(x) = M$. Let us assume that there exits $t_1\geqslant 0$ such that $h(t_1) = 1/(c t_1/2 + 1/M)$ and $t_2>t_1$ such that $h(t) > 1/(c t/2 + 1/M)$ for all $t \in [t_1,t_2]$. Then, for all $t \in [t_1,t_2]$, we have $\chi_\varepsilon(h(t))=0$ and 
$$h'(t)\ \leqslant - \half\, c\, \frac{1}{\left(c\, t/2\, +\, 1/M\right)^2}, \quad \implies \quad h(t)\ \leqslant\ \frac{1}{c\, t/2\, +\, 1/M}. $$
The invertibility of the map $x \mapsto X(t,x)$ ends the proof of \eqref{Oleinikellep}.
\qed

The inequality \eqref{uxenergy} is not uniform on $\ell$ (when $\ell$ is small). In order to obtain a strong limit as $\ell \to 0$ one need to obtain a uniform estimate of the solution. The following lemma gives a uniform estimate of the total variation of $u^{\ell,\varepsilon}$ which plays a very important role in the study of the limiting case $\ell \to 0$ in Section \ref{sec:zero} below.

\begin{lem}\label{lemTV} Let $f$ be a smooth flux such that $0<c \leqslant f''(u) \leqslant C$.
Let also $u_0 \in H^1$ with $u_0' \in L^1$ and $M \eqdef \sup_{x \in \mathds{R}} u_0'(x) < \infty$. Then
\begin{equation}\label{TVep}
\left\| u^{\ell, \varepsilon} \right\|_{L^\infty}\, \leqslant\, \left\| u_x^{\ell, \varepsilon} \right\|_{L^1}\, \leqslant\ \left\| u_0' \right\|_{L^1} \left( c\, M\, t/2\, +\, 1 \right)^{2\, C/c}, \qquad \forall \ell,\varepsilon \in (0,\infty).
\end{equation}
\end{lem}
\proof
Multiplying \eqref{rBep_x} by $\mathrm{sign}(q^{\ell,\varepsilon})$, we obtain
\begin{align*}
\left|q^{\ell,\varepsilon}\right|_{t}\, +\, \left[ f'(u^{\ell,\varepsilon})\, \left|q^{\ell,\varepsilon} \right| \right]_{x}\ 
&=\ - \ell^2\, \mathrm{sign}(q^{\ell,\varepsilon})\, P^\varepsilon_{xx}\\
&=\ - \ell^2\, \mathds{1}_{q^{\ell,\varepsilon} >0} P^\varepsilon_{xx}\ +\ \ell^2\, \mathds{1}_{q^{\ell,\varepsilon} \leqslant 0} P^\varepsilon_{xx}\\
&=\ - 2\, \ell^2\, \mathds{1}_{q^{\ell,\varepsilon} >0} P^\varepsilon_{xx}\ +\ \ell^2\,  P^\varepsilon_{xx}\\
&=\  2\, \mathds{1}_{q^{\ell,\varepsilon} >0} \left[ \half\, f''(u^{\ell,\varepsilon}) \left( q^{\ell,\varepsilon} \right)^2\, +\, \half\, f''(u^{\ell,\varepsilon})\, \chi_\varepsilon(q^{\ell,\varepsilon})\,  -\, P^\varepsilon \right] +\, \ell^2\,  P^\varepsilon_{xx}\\
&\leqslant\  \mathds{1}_{q^{\ell,\varepsilon} >0}\ f''(u^{\ell,\varepsilon}) \left( q^{\ell,\varepsilon} \right)^2\  +\ \ell^2\,  P^\varepsilon_{xx}\\
&\leqslant\  \frac{C}{c\, t/2\, +\, 1/M}\, 
  \left| q^{\ell,\varepsilon} \right|  +\ \ell^2\,  P^\varepsilon_{xx}.
\end{align*}
Integrating with respect to $x$ we obtain 
\begin{equation*}
\frac{\mathrm{d}}{\mathrm{d}t}\, \int_\mathds{R} \left|q^{\ell,\varepsilon}\right| \mathrm{d}x\ \leqslant\ \frac{C}{c\, t/2\, +\, 1/M}\,  \int_\mathds{R} \left|q^{\ell,\varepsilon}\right| \mathrm{d}x.
\end{equation*}
The last inequality with Gronwall lemma imply \eqref{TVep}.
 \qed
 
 In the study the limiting cases $\varepsilon \to 0$ and $\ell \to \infty$ in Sections \ref{sec:precompactness} and \ref{sec:infty} below, we use Young measures to identify the limits of some nonlinear terms depending on the derivative of the solution.  
The $L^2$ estimate \eqref{uxenergy} is uniform on $\varepsilon$ and $\ell$ (when $\ell$ is far from $0$). 
This is not enough to take the limit of the quadratic terms, in that case one need an estimate in a $L^p$ space for some $p>2$.
In the following lemma we present an estimate of the $L^p$ norm of $u^{\ell,\varepsilon}_x$ for any $p<3$.

\begin{lem}\label{lem:alpha+2}
Let $\ell>0$, $\alpha \in (0,1)$, $T>0$ and $ [a,b] \subset \R$, then there exists a constant $C=C(\alpha,T,a,b,\ell)>0$, such that for all $\varepsilon>0$ we have
\begin{equation}\label{alpha+2}
\int_0^T \int_a^b \left[ |u^{\ell,\varepsilon}_t|^{2+\alpha}\ +\ |u^{\ell,\varepsilon}_x|^{2+\alpha} \right] \mathrm{d}x\, \mathrm{d}t\ \leqslant\ C.
\end{equation}
\end{lem}
\begin{rem}
The constant $C>0$ in \eqref{alpha+2} depends on $\ell$. 
In Lemma \ref{lem:alpha+22} below, we prove that if $\ell$ is far from $0$, one can chose a constant $C$ independent on $\ell \geqslant 1$.
\end{rem}
\proof
Without losing generality, we consider that $\alpha = 2k/(2k+1)$ where $k \in \mathds{N}^*$.
Multiplying \eqref{rBep_x} by $\left(q^{\ell,\varepsilon}\right)^\alpha$, we obtain 
\begin{gather} \nonumber
\left[{\textstyle\frac{\left(q^{\ell,\varepsilon}\right)^{\alpha+1}}{\alpha + 1}} \right]_{t} +\
\left[{\textstyle \frac{ f'(u^{\ell,\varepsilon})\left(q^{\ell,\varepsilon}\right)^{\alpha+1}}{\alpha + 1}} \right]_{x}
+\ {\textstyle \frac{\alpha - 1}{2\, (\alpha+1)}}\, f''(u^{\ell,\varepsilon})\, (q^{\ell,\varepsilon})^{\alpha + 2}\ 
+\ (q^{\ell,\varepsilon})^\alpha\, P^\varepsilon\\ \label{q:alpha+1}
 =\ \half\,  f''(u^{\ell,\varepsilon})\, (q^{\ell,\varepsilon})^\alpha\, \chi_\varepsilon(q^{\ell,\varepsilon})\ 
 \geqslant\ 0.
\end{gather}
Let $\Omega = [0,T] \times [a,b]$, $\tilde{\Omega} = [0, T+1) \times (a-1,b+1)$ and let $\varphi \in C^\infty_c(\tilde{\Omega})$ be a non negative function such that $\varphi(t,x)=1$ on $\Omega$.
Multiplying \eqref{q:alpha+1} by $\varphi(t,x)$ and using integration by parts with the energy conservation \eqref{eneequation} we obtain
\begin{align*}
\int_\Omega (q^{\ell,\varepsilon})^{\alpha + 2}\, \mathrm{d}x\, \mathrm{d}t\ 
\leqslant &\ \int_{\tilde{\Omega}}  \varphi(t,x)\, (q^{\ell,\varepsilon})^{\alpha + 2}\, \mathrm{d}x\, \mathrm{d}t\\
\leqslant &\ {\textstyle \frac{2\, (\alpha+1)}{c\, (1-\alpha)}}\, \int_{\tilde{\Omega}}  \varphi(t,x)\, {\textstyle \frac{1-\alpha}{2\, (\alpha+1)}}\, f''(u^{\ell,\varepsilon})\, (q^{\ell,\varepsilon})^{\alpha + 2}\, \mathrm{d}x\, \mathrm{d}t\\
\leqslant &\ {\textstyle \frac{2\, (\alpha+1)}{c\, (1-\alpha)}}\, \int_{\tilde{\Omega}}  \varphi(t,x)\, (q^{\ell,\varepsilon})^\alpha\, P^\varepsilon\, \mathrm{d}x\, \mathrm{d}t\ -{\textstyle \frac{\alpha+1}{c\, (1-\alpha)}}\, \int_{a-1}^{b+1} {\textstyle \frac{(\partial_x u_0^\varepsilon)^{\alpha+1}}{\alpha + 1}}\, \varphi(0,x)\, \ud x
\\
&- {\textstyle \frac{2\, (\alpha+1)}{c\, (1-\alpha)}}\, \int_{\tilde{\Omega}}  {\textstyle\frac{\left(q^{\ell,\varepsilon}\right)^{\alpha+1}}{\alpha + 1}}  \left\{\varphi_t(t,x)\, +\, \varphi_x(t,x)\, f'(u^{\ell,\varepsilon})\  \right\} \mathrm{d}x\, \mathrm{d}t\\
\leqslant &\ C \left[ \|q^{\ell,\varepsilon}\|_{L^\infty_t L^2_x}^\alpha \|P^\varepsilon\|_{L^\infty_t L^{\frac{2}{2-\alpha}}_x}\, +\, 1\, +\, \|q^{\ell,\varepsilon}\|_{L^\infty_t L^2_x}^{\alpha+1} \left( \|f'\! \left(u^{\ell,\varepsilon}\right)\|_{L^\infty(\tilde{\Omega})} +\, 1 \right)   \right].
\end{align*}
Using \eqref{PLp}, \eqref{uLinf} and the energy conservation \eqref{eneequation}, we obtain  
\begin{align*}
\int_\Omega (q^{\ell,\varepsilon})^{\alpha + 2}\, \mathrm{d}x\, \mathrm{d}t\ 
\leqslant \ C.
\end{align*}
Then \eqref{alpha+2} follows directly from \eqref{rBep} and \eqref{PLp}.
\qed

\section{Precompactness of the approximated equation}\label{sec:precompactness}

The aim of this section is to prove Theorem \ref{thm:existence}. For that purpose, we fix $\ell>0$ and we study the limiting case $\varepsilon \to 0$. 
\begin{lem}\label{Strong_conv}
There exist $u^\ell \in L^\infty ([0,\infty), H^1(\mathds{R}))$ and a subsequence of $(u^{\ell,\varepsilon})_\varepsilon$ noted also $(u^{\ell,\varepsilon})_\varepsilon$ such that, as $\varepsilon \to 0$, we have
\begin{align*}
u^{\ell,\varepsilon} \quad &\to \quad u^\ell \qquad \mathrm{in}\  L^{\infty}_{loc}([0,\infty) \times \mathds{R}), \\
u^{\ell,\varepsilon} \quad &\rightharpoonup \quad u^\ell \qquad \mathrm{in}\   H^1([0,T]\times \mathds{R}),\ \forall T>0.
\end{align*}
\end{lem}
\proof Using the energy equation \eqref{eneequation}, we obtain that $u^{\ell,\varepsilon}$ is uniformly (on $\varepsilon$) bounded in $L^\infty([0,\infty), H^1(\mathds{R}))$. Then \eqref{rBep}, \eqref{uLinf} with \eqref{PLp} imply
\begin{equation}\label{L2L2}
\left\|u_t^{\ell,\varepsilon}  \right\|_{L^2([0,T] \times \mathds{R})}\ \leqslant\ C_{T,\ell}.
\end{equation}
Then, the weak convergence in $H^1([0,T] \times \mathds{R})$ follows directly. Using the inequality 
\begin{equation*}
\left\| u^{\ell,\varepsilon}(t,\cdot)\, -\, u^{\ell,\varepsilon}(s,\cdot) \right\|_{L^2(\mathds{R})}^2\ =\ \int_\mathds{R} \left( \int_s^t u^{\ell,\varepsilon}_t(\tau,x)\, \mathrm{d}\tau \right)^2 \mathrm{d}x\ \leqslant\ |t-s|\, \left\|u^{\ell,\varepsilon}_t \right\|_{L^2([0,T] \times \mathds{R})}^2,
\end{equation*}
with \eqref{L2L2} we obtain that for any $\ell > 0$ we have 
\begin{equation*}
\lim_{t \to s} \left\| u^{\ell,\varepsilon}(t,\cdot)\, -\, u^{\ell,\varepsilon} (s,\cdot) \right\|_{L^2(\mathds{R})}\ =\ 0
\end{equation*}
uniformly on $\varepsilon$. Then, using Theorem 5 in \cite{Simon} we can prove that up to a subsequence,
$u^{\ell,\varepsilon}$ converges uniformly to $u^\ell$ on any compact set of $[0,\infty) \times \mathds{R} $ as $\varepsilon \to 0$.
\qed

In order to obtain the precompactness of $(P^\varepsilon)_\varepsilon$, from \eqref{Pepdef}, we write $P^\varepsilon = P_1^\varepsilon + P_2^\varepsilon$ such that 
\begin{equation*}
P_1^\varepsilon\ \eqdef\ \half\,\mathfrak{G}\ast \left\{ f''(u^{\ell,\varepsilon})\, (u^{\ell,\varepsilon}_x)^{\,2} \right\}, \qquad P_2^\varepsilon\ \eqdef\ \half\,\mathfrak{G}\ast \left\{ f''(u^{\ell,\varepsilon})\, \chi_\varepsilon(u_x^{\ell,\varepsilon})  \right\}.
\end{equation*}

\begin{lem}\label{Strong_conv1.5}
There exist $\widetilde{P} \in L^\infty ([0,\infty), H^1(\mathds{R}))$ and a subsequence of $( P^\varepsilon_1)_\varepsilon$ noted also $( P^\varepsilon_1)_\varepsilon$ such that, as $\varepsilon \to 0$, we have
\begin{align}\label{P1conv}
P_1^\varepsilon \quad &\to \quad \widetilde{P} \qquad \mathrm{in}\   L^p_{loc}([0,\infty) \times \mathds{R}),\ \forall p \in (1,\infty).
\end{align}
Moreover, there exists a constant $C_\ell >0$, such that for all $\varepsilon >0$ we have
\begin{equation}\label{P2conv}
\left\|f''(u^{\ell,\varepsilon})\, \chi_\varepsilon (q^{\ell,\varepsilon}) \right\|_{L^1([0,\infty) \times \R)}\, +\, \left\| P_2^\varepsilon \right\|_{L^1([0,\infty), W^{2,1}(\R))}\, \leqslant\ \varepsilon\, C_\ell.
\end{equation}
\end{lem}
\proof
We recall that $\chi_\varepsilon(\xi) = 0 $ for any $\xi \in [-1/\varepsilon, \infty)$ and that $f''(u) \geqslant 0$, $\chi_\varepsilon(\xi) \geqslant 0$ for any $u,\xi \in \R$. Then, using the energy equation \eqref{eneequation} we obtain 
\begin{align*}
\int_{\mathds{R}^+} \int_\mathds{R} f''(u^{\ell,\varepsilon})\, \chi_\varepsilon (q^{\ell,\varepsilon})\, \mathrm{d}x\, \mathrm{d}t\ 
&=\  \int_{\{q^{\ell,\varepsilon} \leqslant -1/\varepsilon\}} f''(u^{\ell,\varepsilon})\, \chi_\varepsilon (q^{\ell,\varepsilon})\, \mathrm{d}x\, \mathrm{d}t\\
&\leqslant \ \int_{\{ q^{\ell,\varepsilon} \leqslant -1/\varepsilon \}}\! \left[  - \varepsilon\, q^{\ell,\varepsilon} \right] \! f''(u^{\ell,\varepsilon})\, \chi_\varepsilon (q^{\ell,\varepsilon})\,  \mathrm{d}x\, \mathrm{d}t \\
&\leqslant \ - \varepsilon \int_{\mathds{R}^+} \int_\mathds{R}  f''(u^{\ell,\varepsilon})\, q^{\ell,\varepsilon}\, \chi_\varepsilon (q^{\ell,\varepsilon})\, \mathrm{d}x\, \mathrm{d}t\ \leqslant\ \varepsilon\, C_\ell.
\end{align*}
This with the definition of $P_2^\varepsilon$ imply \eqref{P2conv}.
From \eqref{rBep} and \eqref{rBep_x} we obtain 
\begin{gather}\nonumber
\left[f''(u^{\ell,\varepsilon})\, \left(q^{\ell,\varepsilon}\right)^2\right]_{t}\, +\  \left[ f''(u^{\ell,\varepsilon})\, f'(u^{\ell,\varepsilon})\, \left(q^{\ell,\varepsilon}\right)^2 \right]_{x}\\
 =\ -\ \ell^2\, f'''(u^{\ell,\varepsilon})\, P_x^\varepsilon \left(q^{\ell,\varepsilon}\right)^2\, -\ 2\, f''(u^{\ell,\varepsilon})\, P^\varepsilon\, q^{\ell,\varepsilon}\
+\ 
 \chi_\varepsilon(q^{\ell,\varepsilon})\, f''(u^{\ell,\varepsilon})^2\, q^{\ell,\varepsilon}.
\end{gather}
Then from the definition of $P^\varepsilon$, the energy equation \eqref{eneequation}, \eqref{uLinf} and \eqref{PLp}  we obtain that the quantity $\left[f''(u^{\ell,\varepsilon}) \left(q^{\ell,\varepsilon}\right)^2\right]_{t}$ is bounded in $L^1([0, \infty), W^{-1,1}_{loc}(\R))$ for any fixed $\ell >0$.
Then, $\partial_t P^\varepsilon_1$ is bounded in $L^1([0, \infty), W_{loc}^{1,1}(\R))$. 
As in \eqref{PLp}, one can prove that $P_1^\varepsilon$ is bounded in $L^\infty([0,\infty), W^{1,\infty}(\R))$.
Using the compact embedding $W_{loc}^{1,\infty}(\R) \Subset L^p_{loc}(\R)$, the continuous embedding $L^p_{loc}(\R) \hookrightarrow L^1_{loc}(\R)$ and Aubin lemma, one obtains the convergence \eqref{P1conv}.
\qed

We proved in Lemma \ref{Strong_conv} the strong precompactness of $(u^{\ell,\varepsilon})_\varepsilon$. However, the best estimation obtained for $q^{\ell,\varepsilon}=u_x^{\ell,\varepsilon}$ is \eqref{alpha+2} which ensures only the weak convergence in $L^p$ with $p< 3$. This is not enough to identify the limit of the quadratic term $(q^{\ell,\varepsilon})^2$.
In order to overcome this problematic, we use the Young measures and we show later that the Young measure obtained at the limit is nothing but a Dirac mass. This insures that the limit of the quadratic term $(q^{\ell,\varepsilon})^2$ is the square of the weak limit of $q^{\ell,\varepsilon}$.

\begin{lem}\label{lemYoungm}
There exist a subsequence of $\left\{q^{\ell,\varepsilon} \right\}_\varepsilon$ denoted also $\left\{q^{\ell,\varepsilon} \right\}_\varepsilon$ and a family of probability Young measures $\mu_{t,x}^\ell$ on $\mathds{R}$, such that for all functions $g \in C(\mathds{R})$ with $g(\xi) = \mathcal{O} (|\xi|^{2})$ at infinity, and for all $\varphi \in C^\infty_c((0,\infty) \times \mathds{R})$ we have
\begin{equation}\label{convegrgenceg}
\lim_{\varepsilon \to 0}\, \int_{\mathds{R}^+ \times \mathds{R}} \varphi(t,x)\, g(q^{\ell,\varepsilon})\, \mathrm{d}x\, \mathrm{d}t\ 
=\ \int_{\mathds{R}^+ \times \mathds{R}} \varphi(t,x) \int_{\mathds{R}} g(\xi)\,  \mathrm{d} \mu_{t,x}^\ell(\xi)\, \mathrm{d}x\, \mathrm{d}t.
\end{equation}
Moreover, the map
\begin{equation}\label{Fatou1}
(t,x)\ \mapsto\  \int_\mathds{R} \xi^2\, \mathrm{d} \mu_{t,x}^\ell (\xi)
\end{equation}
belongs to $L^\infty(\mathds{R}^+, L^1(\mathds{R})) $.
\end{lem}
\proof If $g(\xi)=\smallO(|\xi|^2)$, then the result is a direct consequences of the energy equation \eqref{eneequation} and  Lemma \ref{lem:Young}.
If $g(\xi)=\mathcal{O}(|\xi|^2)$, let $\psi$ be a smooth cut-off function with $\psi(\xi)=1$ for $|\xi| \leqslant 1$ and $\psi(\xi)=0$ for $|\xi| \geqslant 2$, then 
\begin{equation}\label{convegrgencegpsi}
\lim_{\varepsilon \to 0}\, \int_{\mathds{R}^+ \times \mathds{R}} \varphi(t,x)\, g_k(q^{\ell,\varepsilon})\,  \mathrm{d}x\, \mathrm{d}t\ 
=\ \int_{\mathds{R}^+ \times \mathds{R}} \varphi(t,x) \int_{\mathds{R}} g_\kappa(\xi)\,  \mathrm{d} \mu_{t,x}^\ell(\xi)\, \mathrm{d}x\, \mathrm{d}t,
\end{equation}
where $g_\kappa(\xi) \eqdef g(\xi) \psi\left( {\textstyle \frac{\xi}{\kappa}}\right)$ with $\kappa>0$. Using Holder inequality, Lemma \ref{lem:alpha+2} with $\Omega = \mathrm{supp}(\varphi)$ we obtain
\begin{gather*}
\left| \int_{\mathds{R}^+ \times \mathds{R}} \varphi(t,x)\, \left( g(q^{\ell,\varepsilon})\, -\, g_\kappa(q^{\ell,\varepsilon}) \right)\, \mathrm{d}x\, \mathrm{d}t \right|\
\leqslant\ 
\int_{\mathrm{supp}(\varphi) \cap \{|q^{\ell,\varepsilon}| \geqslant \kappa\}} |\varphi(t,x)|\, \left|g(q^{\ell,\varepsilon})\right|\, \mathrm{d}x\, \mathrm{d}t \\
\leqslant\ C \left( \int_{\mathrm{supp}(\varphi)} \left|g(q^{\ell,\varepsilon})\right|^{p/2}\, \mathrm{d}x\, \mathrm{d}t \right)^{2/p} \left( \int_{\mathrm{supp}(\varphi) \cap \{|q^{\ell,\varepsilon}| \geqslant \kappa\}} \mathrm{d}x\, \mathrm{d}t \right)^\frac{p-2}{p}
\\
\leqslant\ C\, \Big[ \big|\! \left\{(t,x) \in \mathrm{supp}(\varphi), |q^{\ell,\varepsilon}| \geqslant \kappa \right\}\! \big| \Big]^\frac{p-2}{p}\
\leqslant\ C\, \kappa^{2-p}.
\end{gather*}
where $2<p<3$.
The last inequality with \eqref{convegrgencegpsi} imply that we can interchange the limits $\kappa \to \infty$ and $\varepsilon \to 0$. Using that $|g_\kappa| \leqslant |g|$ and the dominated convergence theorem we obtain \eqref{convegrgenceg}.
\qed

Now we define
\begin{equation}\label{bardef}
\overline{g(q)}\ \eqdef\ \int_{\mathds{R}} g(\xi)\, \mathrm{d} \mu_{t,x}^\ell(\xi)
\end{equation}
which is from \eqref{convegrgenceg} the weak limit of $g(q^{\ell,\varepsilon})$. 
Our aim is to show that $\overline{q^2}=\overline{q}^2$ which means that $\mu^\ell_{t,x}= \delta_{\overline{q}(t,x)}$.
For that purpose, we need the following technical lemmas.

\begin{lem} As $t \to 0$ we have
\begin{gather} \label{tauto0}
\left\|u^\ell-u_0 \right\|_{H^1(\mathds{R})}\ \to\ 0 \quad \qquad \mathrm{and} \qquad \quad
\int_\mathds{R} \left( \overline{q^2}\ -\ \overline{q}^2\, \right) \mathrm{d}x\ \to\ 0.
\end{gather}
\end{lem}

\proof

From Lemma \ref{Strong_conv} and Lemma \ref{lemYoungm} we have for all $t > 0$
\begin{gather*}
\left(u^{\ell,\varepsilon}(t,\cdot), \ell\, u^{\ell,\varepsilon}_x(t,\cdot) \right)\ \rightharpoonup\ \left(u^\ell(t,\cdot), \ell\, u^\ell_x(t,\cdot) \right) \qquad \mathrm{as}\ \varepsilon\, \to\, 0 \qquad \mathrm{in\ } L^2(\mathds{R}),\\
u_x^{\ell,\varepsilon}(t,\cdot)^2\ \rightharpoonup\ \overline{q^2}(t,\cdot) \qquad \mathrm{when}\ \varepsilon\, \to\, 0 \qquad \mathrm{in\ } \mathcal{D}'(\mathds{R}).
\end{gather*}
This, with Jensen's inequality and the energy equation \eqref{eneequation} imply that 
\begin{align}\nonumber
\left\|\left(u^\ell(t), \ell\, u^\ell_x(t)\right) \right\|_{L^2(\mathds{R})}^2\ &\leqslant\ \left\|\left(u^\ell(t), \ell\, \sqrt{\overline{q^2}}(t)\right) \right\|_{L^2(\mathds{R})}^2\  \leqslant\
\liminf_{\varepsilon \to 0} \left\|\left(u^{\ell,\varepsilon}(t), \ell\, u^{\ell,\varepsilon}_x(t) \right)\right\|_{L^2(\mathds{R})}^2\\ \label{Energyfin}
& \leqslant\ \lim_{\varepsilon \to 0} \left\|\left(u^\varepsilon_0, \ell\, \partial_x u^\varepsilon_0 \right)\right\|_{L^2(\mathds{R})}^2\ =\ \left\|\left(u_0, \ell\,  u'_0 \right)\right\|_{L^2(\mathds{R})}^2.
\end{align}
The energy inequality \eqref{eneequation} implies that the sequences $(u^{\ell,\varepsilon})_\varepsilon$ and $(q^{\ell,\varepsilon})_\varepsilon$ are bounded in the space $L^\infty([0,\infty), L^2(\mathds{R}))$.
Using \eqref{rBep}, we can prove that for all $T>0$ and for all $\varphi \in H^1(\mathds{R})$, the map 
\begin{equation}
t\ \mapsto\ \int_\mathds{R} \varphi(x)\, \left(u^{\ell,\varepsilon}, \ell\, q^{\ell,\varepsilon}\right) \mathrm{d}x
\end{equation}
is uniformly (on $t \in [0,T]$ and $\varepsilon >0 $) continuous. Then, Lemma \ref{lem:C_w} implies that 
\begin{equation}\label{Wconv2}
\left(u^\ell(t,\cdot), \ell\, q^\ell (t,\cdot)\right)\ \rightharpoonup\ \left(u_0, \ell\, u_0'\right)\qquad \mathrm{when}\ t\, \to\, 0 \qquad \mathrm{in\ } L^2(\mathds{R}),
\end{equation}
which implies that 
\begin{equation*}
\left\|\left(u_0, \ell\, u_0'\right) \right\|_{L^2}^2\ \leqslant\ \liminf_{t \to 0}  \left\|\left(u^\ell(t), \ell\, u^\ell_x(t)\right)\right\|_{L^2}^2.
\end{equation*}
On another hand, \eqref{Energyfin} implies
\begin{equation*}
\limsup_{t \to 0}\left\|\left(u^\ell(t), \ell\, u^\ell_x(t)\right) \right\|_{L^2(\mathds{R})}^2\ \leqslant\ \left\|\left(u_0, \ell\,  u'_0 \right)\right\|_{L^2(\mathds{R})}^2,
\end{equation*}
then 
\begin{equation}\label{Normconv}
\lim_{t \to 0}\left\|\left(u^\ell(t), \ell\, u^\ell_x(t)\right) \right\|_{L^2(\mathds{R})}^2\ =\ \left\|\left(u_0, \ell\,  u'_0 \right)\right\|_{L^2(\mathds{R})}^2,
\end{equation}
which implies with \eqref{Wconv2} that 
\begin{equation}\label{Sconv}
\left(u^\ell(t,\cdot), \ell\, q^\ell (t,\cdot)\right)\ \to\ \left(u_0, \ell\, u_0'\right)\qquad \mathrm{when}\ t\, \to\, 0 \qquad \mathrm{in\ } L^2(\mathds{R}).
\end{equation}
The inequality \eqref{Energyfin} with \eqref{Normconv} imply 
\begin{equation}\label{Normconv2}
\lim_{t \to 0}\left\|\left(u^\ell(t), \ell\, \sqrt{\overline{q^2}} (t)\right) \right\|_{L^2(\mathds{R})}^2\ =\
\lim_{t \to 0}\left\|\left(u^\ell(t), \ell\, u^\ell_x(t)\right) \right\|_{L^2(\mathds{R})}^2\ =\ \left\|\left(u_0, \ell\,  u'_0 \right)\right\|_{L^2(\mathds{R})}^2.
\end{equation}
Then \eqref{tauto0} follows directly from \eqref{Sconv} and \eqref{Normconv2}. 

\qed

For any $\kappa >0$, we define
\begin{equation}\label{Sdef}
S_\kappa(\xi)\ \eqdef\ \half\, \xi^2\ -\ \half\, (\xi\, +\, \kappa)^2\, \mathds{1}_{\xi\leqslant-\kappa}\ -\ \half\, (\xi\, -\, \kappa)^2\, \mathds{1}_{\xi\geqslant \kappa}\ =\ 
\begin{cases}
-\kappa\, (\xi\, +\, \half\, \kappa), & \xi \leqslant -\kappa, \\
\half\, \xi^2, & |\xi| \leqslant \kappa, \\
\kappa\, (\xi\, -\, \half\, \kappa), & \xi \geqslant \kappa.
\end{cases}
\end{equation}
\begin{equation}\label{Tdef}
T_\kappa(\xi)\ \eqdef\ S_\kappa'(\xi)\ =\ \xi\ -\ (\xi\, +\, \kappa)\, \mathds{1}_{\xi\leqslant-\kappa}\ -\ (\xi\, -\, \kappa)\, \mathds{1}_{\xi\geqslant \kappa}\ =\ 
\begin{cases}
-\kappa, & \xi \leqslant -\kappa, \\
\xi, & |\xi| \leqslant \kappa, \\
\kappa, & \xi \geqslant \kappa.
\end{cases}
\end{equation}

\begin{lem}\label{lem:ST_k}
For all $T>0$, we have
\begin{gather}\label{L1conv}
\lim_{\kappa \to \infty} \left\| \overline{T_\kappa(q)}\ -\ T_\kappa\left(\overline{q}\right) \right\|_{L^1([0,T] \times \mathds{R})}\ =\
\lim_{\kappa \to \infty} \left\| \overline{T_\kappa(q)}\ -\ \overline{q} \right\|_{L^1([0,T] \times \mathds{R})}\ =\ 0. 
\end{gather}
Moreover, for all $\kappa >0$ we have
\begin{align}\label{TS}
\half\, \left( \overline{T_\kappa(q)}\ -\ T_\kappa\left(\overline{q}\right) \right)^2\ 
&\leqslant\ \overline{S_\kappa(q)}\ -\ S_\kappa\left(\overline{q}\right).
\end{align}
\end{lem}

\proof 
From \eqref{Tdef} we have 
\begin{align*}
\left| T_\kappa\left(\xi\right)\ - \xi\right|\ \leqslant\
\left|\xi\, +\, \kappa \right|\, \mathds{1}_{\xi\leqslant-\kappa}\ +\ \left|\xi\, -\, \kappa \right|\, \mathds{1}_{\xi\geqslant \kappa}\
\leqslant\ 2\, \left|\xi\right|\, \mathds{1}_{\kappa \leqslant \left|\xi\right|}\ 
\leqslant\ {\textstyle \frac{2}{\kappa}}\, \xi^2.
\end{align*}
Then, we have 
\begin{equation}
\left| \overline{T_\kappa\left(q\right)}\ - T_\kappa\left(\overline{q}\right) \right|\ \leqslant\ 
\left| \overline{T_\kappa\left(q\right)}\ - \overline{q} \right|\ +\
\left| T_\kappa\left(\overline{q}\right)\ - \overline{q} \right|\
\leqslant {\textstyle \frac{2}{\kappa}}\, \left( \overline{q^2}\, +\, \overline{q}^2 \right).
\end{equation}
Jenson's inequality imply that $\overline{q}^2 \leqslant \overline{q^2}$. Lemma \ref{lemYoungm} implies that $\overline{q^2} \in L^\infty(\mathds{R}^+ , L^1(\mathds{R}))$. Then \eqref{L1conv} follows directly.

Cauchy--Schwarz inequality implies that $\overline{T_\kappa(q)}^2 \leqslant \overline{T_\kappa(q)^2}$, then, using the definition \eqref{Tdef} we obtain
\begin{align} \nonumber
\left( \overline{T_\kappa(q)}\ -\ T_\kappa\left( \overline{q}\right) \right)^2\ 
\leqslant &\ \overline{T_\kappa(q)^2}\ +\ T_\kappa\left( \overline{q}\right)^2\ -\ 2\, T_\kappa\left( \overline{q}\right)\, \overline{T_\kappa(q)}\\ \nonumber
=&\
\overline{T_\kappa(q)^2}\ +\ T_\kappa\left( \overline{q}\right)^2\ -\ 2\, T_\kappa\left( \overline{q}\right)\, \overline{q}\ +\ 2\, T_\kappa\left( \overline{q}\right)\, \overline{ \left( q\, +\, \kappa \right)\, \mathds{1}_{q \leqslant -\kappa}}\\ \nonumber
&+\ 2\, T_\kappa(\overline{q})\, \overline{ \left( q\, -\, \kappa \right)\, \mathds{1}_{q \geqslant \kappa}}\\ \nonumber
=&\ \overline{T_\kappa(q)^2}\ +\ 2\, T_\kappa\left( \overline{q}\right)\, \left[ \overline{ \left( q\, +\, \kappa \right)\, \mathds{1}_{q \leqslant -\kappa}}\ -\ \left( \overline{q}\, +\, \kappa \right)\, \mathds{1}_{\overline{q} \leqslant - \kappa} \right]\\ \nonumber
& - T_\kappa\left( \overline{q}\right)^2\ +\ 2\, T_\kappa\left( \overline{q}\right)\, \left[ \overline{ \left( q\, -\, \kappa \right)\, \mathds{1}_{q \geqslant \kappa}}\ -\ \left( \overline{q}\, -\, \kappa \right)\, \mathds{1}_{\overline{q} \geqslant \kappa} \right]\\ \nonumber
\leqslant &\ \overline{T_\kappa(q)^2}\ -\ 2\, \kappa\, \left[ \overline{ \left( q\, +\, \kappa \right)\, \mathds{1}_{q \leqslant -\kappa}}\ -\ \left( \overline{q}\, +\, \kappa \right)\, \mathds{1}_{\overline{q} \leqslant - \kappa} \right]\\ \label{TT} 
& - T_\kappa\left( \overline{q}\right)^2\ +\ 2\, \kappa\, \left[ \overline{ \left( q\, -\, \kappa \right)\, \mathds{1}_{q \geqslant \kappa}}\ -\ \left( \overline{q}\, -\, \kappa \right)\, \mathds{1}_{\overline{q} \geqslant \kappa} \right],
\end{align}
where the last inequality follows from Jensen's inequality with the concavity of $\xi \mapsto \left( \xi + \kappa \right) \mathds{1}_{\xi \leqslant - \kappa}$, the convexity of $\xi \mapsto \left( \xi - \kappa \right) \mathds{1}_{\xi \geqslant  \kappa}$ and $-\kappa \leqslant T_\kappa(\xi)\leqslant \kappa$. Since 
\begin{equation*}
S_\kappa(\xi)\ =\ \half\, T_\kappa(\xi)^2\ +\ \kappa\, \left( \xi\, -\, \kappa \right)\, \mathds{1}_{\xi \geqslant  \kappa}\ -\ \kappa\, \left( \xi\, +\, \kappa \right)\, \mathds{1}_{\xi \leqslant - \kappa}
\end{equation*}
we have 
\begin{align*}
\overline{S_\kappa(q)}\ &=\ \half\, \overline{T_\kappa(q)^2}\ +\ \kappa\, \overline{ \left( q\, -\, \kappa \right)\, \mathds{1}_{q \geqslant  \kappa}}\ -\ \kappa\, \overline{ \left( q\, +\, \kappa \right)\, \mathds{1}_{q \leqslant - \kappa}},\\
S_\kappa\left(\overline{q}\right)\ &=\ \half\, T_\kappa(\overline{q})^2\ +\ \kappa\, \left( \overline{q}\, -\, \kappa \right)\, \mathds{1}_{\overline{q} \geqslant  \kappa}\ -\ \kappa\, \left( \overline{q}\, +\, \kappa \right)\, \mathds{1}_{\overline{q} \leqslant - \kappa}.
\end{align*}
The last two identities with \eqref{TT} imply \eqref{TS}.
\qed

Now we are ready to prove the main lemma of this section that is the key to pass to the limit $\varepsilon \to 0$ in the quadratic terms of \eqref{rBep}.

\begin{lem}\label{lem:2Dirac}
 The measure $\mu_{t,x}^\ell$ given in Lemma \ref{lemYoungm} is a Dirac measure, and
\begin{equation}
\mu_{t,x}^\ell(\xi)\ =\ \delta_{u_x^\ell(t,x)}(\xi).
\end{equation}
\end{lem}
\proof 
The Idea of the proof is to consider the defect measure $\Delta = \overline{q^2} - \overline{q}^2 \geqslant 0$, we obtain then a transport inequation for $\sqrt{\Delta}$ (see \eqref{Tg} below). Then, since $\Delta(t=0)=0$ we deduce that $\Delta = 0$ for all time.

However, when multiplying \eqref{rBep_x} by $q^{\ell,\varepsilon}$ a cubic term on $q^{\ell,\varepsilon}$ appears, in that case Lemma \ref{lemYoungm} cannot be used. That is the reason of using the cut-off functions defined in \eqref{Sdef} and \eqref{Tdef} and then we take the limit $\kappa \to \infty$. The proof is given in several steps.

\textbf{Step 1.} Multiplying \eqref{rBep_x} by $T_\kappa(q^{\ell,\varepsilon})$ we obtain 
\begin{gather}\nonumber
\left[S_\kappa \left(q^{\ell,\varepsilon}\right)\right]_{t}\, +\, \left[ f'(u^{\ell,\varepsilon})\, S_\kappa \left( q^{\ell,\varepsilon}\right) \right]_{x}\, =\ -\ P^\varepsilon\, T_\kappa\left(q^{\ell,\varepsilon}\right)\\
-\ \half
\left[  (q^{\ell,\varepsilon})^2\ -\  \chi_\varepsilon(q^{\ell,\varepsilon}) \right] f''(u^{\ell,\varepsilon})\, T_\kappa\left(q^{\ell,\varepsilon}\right)\ +\ q^{\ell,\varepsilon}\, f''(u^{\ell,\varepsilon})\, S_\kappa\left(q^{\ell,\varepsilon}\right).
\end{gather}
Using Lemma \ref{Strong_conv1.5} and taking $\varepsilon \to 0$ we obtain
\begin{equation}\label{Sbar0}
\left[ \overline{S_\kappa \left(q\right)}\right]_{t}\, +\, \left[ f'(u^\ell)\, \overline{S_\kappa \left( q\right)} \right]_{x}\, =\ -\ \overline{T_\kappa\left(q\right)}\, \widetilde{P}\
+\ \half\, f''(u^\ell) \left[2\, \overline{q\, S_\kappa\left(q\right)}\, -\, \overline{q^2\, T_\kappa\left(q\right)}  \right].
\end{equation}

\textbf{Step 2.} Taking $\varepsilon \to 0$ in \eqref{rBep_x} we obtain
\begin{equation}\label{barqeq}
\overline{q}_{t}\, +\, \left[f'(u^\ell)\, \overline{q}\right]_{x}\  =\ -\ \widetilde{P}\, +\ \half\, f''(u^\ell)\, \overline{q^2}.
\end{equation}
Let $j_\varepsilon$ be a Friedrichs mollifier and $\overline{q}^\varepsilon \eqdef \overline{q} \ast j_\varepsilon$ then using Lemma \ref{lemRenor} we obtain
\begin{equation*}
\overline{q}^\varepsilon_{t}\, +\, \left[f'(u^\ell)\, \overline{q}^\varepsilon\right]_{x}\  =\ -\ \widetilde{P} \ast j_\varepsilon\, +\ \half\, f''(u^\ell)\, \overline{q^2}\ +\ \theta_\varepsilon,
\end{equation*}
where $\theta_\varepsilon \to 0$ in $L^1_{loc}([0,\infty) \times \R)$ as $\varepsilon \to 0$. Multiplying by $T_\kappa\left(\overline{q}^\varepsilon\right)$ we obtain
\begin{gather*}
\left[S_\kappa\left(\overline{q}^\varepsilon\right)\right]_{t}\, +\, \left[f'(u^\ell)\, S_\kappa\left(\overline{q}^\varepsilon\right)\right]_{x}\ =\ - \left\{\widetilde{P} \ast j_\varepsilon \right\} T_\kappa\left(\overline{q}^\varepsilon\right) \ +\ \theta_\varepsilon\, T_\kappa\left(\overline{q}^\varepsilon\right)\\
+\ f''(u^\ell)\, \overline{q}\, S_\kappa\left(\overline{q}^\varepsilon\right)\, +\ \half\, f''(u^\ell)\, \overline{q^2}\, T_\kappa\left(\overline{q}^\varepsilon\right)\, -\ f''(u^\ell)\, \overline{q}\, \overline{q}^\varepsilon\, T_\kappa\left(\overline{q}^\varepsilon\right).
\end{gather*}
Taking $\varepsilon \to 0$ we obtain 
\begin{gather}\nonumber
\left[S_\kappa\left(\overline{q}\right)\right]_{t}\, +\, \left[f'(u^\ell)\, S_\kappa\left(\overline{q}\right)\right]_{x}\ =\ -\, \widetilde{P}\, T_\kappa\left(\overline{q}\right)\\ \label{Sbar1}
+\ f''(u^\ell)\, \overline{q}\, S_\kappa\left(\overline{q}\right)\, -\ \half\, f''(u^\ell)\, \overline{q}^2\, T_\kappa\left(\overline{q}\right)\, +\ \half\, f''(u^\ell)\left( \overline{q^2}\, -\, \overline{q}^2 \right) T_\kappa\left(\overline{q}\right).
\end{gather}

\textbf{Step 3.} From \eqref{Sbar0} and \eqref{Sbar1} we obtain
\begin{gather}\nonumber
\left[ \overline{S_\kappa \left(q\right)}\, -\, S_\kappa\left(\overline{q}\right)\right]_{t}\, +\, \left[ f'(u^\ell) \left( \overline{S_\kappa \left( q\right)}\, -\, S_\kappa\left(\overline{q}\right) \right) \right]_{x}\, =\ -\,  \widetilde{P} \left( \overline{T_\kappa\left(q\right)}\, -\, T_\kappa\left(\overline{q}\right) \right)\\ \label{Sbar2}
+\ \half\, f''(u^\ell) \left[2\, \overline{q\, S_\kappa\left(q\right)}\, -\, \overline{q^2\, T_\kappa\left(q\right)}\ -\ 2\, \overline{q}\, S_\kappa\left(\overline{q}\right)\, 
+\ \overline{q}^2\, T_\kappa\left(\overline{q}\right)\, +\, \left(\overline{q}^2\, -\, \overline{q^2} \right) T_\kappa\left(\overline{q}\right)\right].
\end{gather}
From \eqref{Sdef} and \eqref{Tdef} we have
\begin{align}\nonumber
\xi^2\, T_\kappa(\xi)\ -\ 2\, \xi\, S_\kappa(\xi)\ 
=&\ 
\xi^2\, T_\kappa(\xi)\ -\ 2\, \xi\, S_\kappa(\xi)\ +\xi^3\ -\xi^3\\ \nonumber
=&\ \xi^2\, \left[T_\kappa(\xi)\, -\, \xi \right]\ +\ \xi\, (\xi\, +\, \kappa)^2\, \mathds{1}_{\xi \leqslant -\kappa}\ +\ \xi\, (\xi-\kappa)^2\, \mathds{1}_{\xi \geqslant \kappa}\\ \nonumber
=&\ \kappa^2\, \left[T_\kappa(\xi)\, -\, \xi \right]\ -\ \left( \xi^2\, -\, \kappa^2\right) \left[ (\xi\, +\, \kappa)\, \mathds{1}_{\xi \leqslant -\kappa}\ +\ (\xi-\kappa)\, \mathds{1}_{\xi \geqslant \kappa} \right]\\ \nonumber
&+\ \xi\, (\xi\, +\, \kappa)^2\, \mathds{1}_{\xi \leqslant -\kappa}\ +\ \xi\, (\xi-\kappa)^2\, \mathds{1}_{\xi \geqslant \kappa}\\ \label{xiTS}
=&\ \kappa^2\, \left[T_\kappa(\xi)\, -\, \xi \right]\ +\ \kappa\, (\xi\, +\, \kappa)^2\, \mathds{1}_{\xi \leqslant -\kappa}\ -\ \kappa\, (\xi-\kappa)^2\, \mathds{1}_{\xi \geqslant \kappa}.
\end{align}
Then from \eqref{Sdef} we have 
\begin{gather}\nonumber
2\, \overline{q\, S_\kappa\left(q\right)}\, -\, \overline{q^2\, T_\kappa\left(q\right)}\ -\ 2\, \overline{q}\, S_\kappa\left(\overline{q}\right)\, 
+\ \overline{q}^2\, T_\kappa\left(\overline{q}\right)\, +\, \left(\overline{q}^2\, -\, \overline{q^2} \right) T_\kappa\left(\overline{q}\right)
\\ \nonumber
=\ 
 \left( T_\kappa\left(\overline{q}\right)\, +\, \kappa \right)\, (\overline{q}\, +\, \kappa)^2\, \mathds{1}_{\overline{q} \leqslant -\kappa}\
+\ \left( T_\kappa\left(\overline{q}\right)\, -\, \kappa \right)\, (\overline{q}\, -\, \kappa)^2\, \mathds{1}_{\overline{q} \geqslant \kappa}  \\ \nonumber
 -\left( T_\kappa\left(\overline{q}\right)\, +\, \kappa \right)\, \overline{(q\, +\, \kappa)^2\, \mathds{1}_{q \leqslant -\kappa}}\
-\ \left( T_\kappa\left(\overline{q}\right)\, -\, \kappa \right)\, \overline{(q\, -\, \kappa)^2\, \mathds{1}_{q \geqslant \kappa}}\\ \label{Step3_1}
-\ \kappa^2\, \left( \overline{T_\kappa(q)}\, -\, T_\kappa\left(\overline{q}\right) \right)\ 
-\ 2\, T_\kappa \left(\overline{q} \right) \left( \overline{S_\kappa(q)}\, -\, S_\kappa\left(\overline{q}\right) \right).
\end{gather} 
From the definition \eqref{Tdef} we have 
\begin{equation}\label{Step3_15}
\left( T_\kappa\left(\overline{q}\right)\, +\, \kappa \right)\, (\overline{q}\, +\, \kappa)^2\, \mathds{1}_{\overline{q} \leqslant -\kappa}\
=\ \left( T_\kappa\left(\overline{q}\right)\, -\, \kappa \right)\, (\overline{q}\, -\, \kappa)^2\, \mathds{1}_{\overline{q} \geqslant \kappa}\ =\ 0.
\end{equation}
Since $T_\kappa\left(\overline{q}\right) \geqslant - \kappa $, then 
\begin{equation}\label{Step3_2}
 -\left( T_\kappa\left(\overline{q}\right)\, +\, \kappa \right)\, \overline{(q\, +\, \kappa)^2\, \mathds{1}_{q \leqslant -\kappa}}\ \leqslant\ 0.
\end{equation}
Let $t_0>0$ and $\kappa \geqslant 2/(c t_0)$, then from Lemma \ref{Lem:Oleinik}, we have for all $t \geqslant t_0$ that $q^{\ell,\varepsilon} \leqslant \kappa $ and $\overline{q} \leqslant \kappa$. Then, using the convexity of $T_\kappa$ on $(-\infty,\kappa)$ and the Jensen's inequality we obtain 
\begin{equation}\label{Step3_3}
-\kappa^2\, \left( \overline{T_\kappa(q)}\, -\, T_\kappa\left(\overline{q}\right) \right)\ \leqslant\ 0, \qquad \forall t \geqslant t_0, \quad   \kappa \geqslant 2/(c\, t_0).
\end{equation}
We take again $t_0>0$ and $\kappa \geqslant 2/(c t_0)$, then for all $\varphi \in C^\infty_c((t_0,\infty) \times \mathds{R})$ we have 
\begin{equation}\label{Step3_4}
\int \overline{(q\, -\, \kappa)^2\, \mathds{1}_{q \geqslant \kappa}}\, \varphi\, \mathrm{d}x\, \mathrm{d}t\ =\ \lim_{\varepsilon \to 0} \int (q^{\ell,\varepsilon}\, -\, \kappa)^2\, \mathds{1}_{q^{\ell,\varepsilon} \geqslant \kappa}\, \varphi\, \mathrm{d}x\, \mathrm{d}t\ =\ 0.
\end{equation}
Defining $\Delta_\kappa \eqdef \overline{S_\kappa \left(q\right)}\, -\, S_\kappa\left(\overline{q}\right) $ and summing up \eqref{Sbar2}, \eqref{Step3_1}, \eqref{Step3_15}, \eqref{Step3_2}, \eqref{Step3_3} and \eqref{Step3_4} we obtain that $\forall t_0>0$, $\forall t \geqslant t_0$ and $\forall \kappa \geqslant 2/(c t_0)$, we have
\begin{equation*}
\left[\Delta_\kappa\right]_t\ +\ \left[ f'(u^\ell)\, \Delta_\kappa \right]_x\ \leqslant\ -\,  \widetilde{P} \left( \overline{T_\kappa\left(q\right)}\, -\, T_\kappa\left(\overline{q}\right) \right)\, -\ f''(u^\ell)\, T_\kappa \left(\overline{q} \right) \Delta_\kappa .
\end{equation*}

\textbf{Step 4.}
Defining $\Delta^\varepsilon_\kappa \eqdef \Delta_\kappa \ast j_\varepsilon$ and using Lemma \ref{lemRenor} we obtain 
\begin{equation*}
\left[\Delta^\varepsilon_\kappa\right]_t\ +\ \left[ f'(u^\ell)\, \Delta^\varepsilon_\kappa \right]_x\ \leqslant\ -\, \widetilde{P} \left( \overline{T_\kappa\left(q\right)}\, -\, T_\kappa\left(\overline{q}\right) \right)\, -\ f''(u^\ell)\, T_\kappa \left(\overline{q} \right) \Delta_\kappa^\varepsilon\ +\ \tilde{\theta}_\varepsilon,
\end{equation*}
where $\tilde{\theta}_\varepsilon \to 0$ as $\varepsilon \to 0$ in $L^1_{loc}((0,\infty)\times \R)$.
Let $\beta >0$, multiplying by $\left( \Delta^\varepsilon_\kappa + \beta \right)^{-1/2}/2 $ we obtain 
\begin{gather*}\nonumber
\left[\sqrt{\Delta^\varepsilon_\kappa + \beta} \right]_t\ +\ \left[ f'(u^\ell)\, \sqrt{\Delta^\varepsilon_\kappa + \beta} \right]_x\ \leqslant\ \half\, f''(u^\ell)\, \frac{\overline{q}\, -\, T_\kappa \left( \overline{q} \right)}{\sqrt{\Delta^\varepsilon_\kappa + \beta}}\, \Delta^\varepsilon_\kappa\ +\ \frac{2\, \beta\, \overline{q}\, f''(u^\ell)\, +\, \tilde{\theta}_\varepsilon}{2\, \sqrt{\Delta^\varepsilon_\kappa + \beta}}\\
 -\, \half\, \widetilde{P}\, \frac{\overline{T_\kappa\left(q\right)}\, -\, T_\kappa\left(\overline{q}\right)}{\sqrt{\Delta^\varepsilon_\kappa + \beta}}.
\end{gather*}
Taking $\varepsilon \to 0$ we obtain 
\begin{gather}\nonumber
\left[\sqrt{\Delta_\kappa + \beta} \right]_t\ +\ \left[ f'(u^\ell)\, \sqrt{\Delta_\kappa + \beta} \right]_x\ \leqslant\ \half\, f''(u^\ell)\, \frac{\overline{q}\, -\, T_\kappa \left( \overline{q} \right)}{\sqrt{\Delta_\kappa + \beta}}\, \Delta_\kappa\ +\ \frac{\beta\, \overline{q}\, f''(u^\ell)}{ \sqrt{\Delta_\kappa + \beta}}\\ \label{Delta_k}
 -\, \half\, \widetilde{P}\, \frac{\overline{T_\kappa\left(q\right)}\, -\, T_\kappa\left(\overline{q}\right)}{\sqrt{\Delta_\kappa + \beta}}.
\end{gather}
Using that $|T_\kappa(\xi)| \leqslant |\xi|$ and $|S_\kappa(\xi)| \leqslant \xi^2/2$ we obtain 
\begin{equation*}
\left| \half\, f''(u^\ell)\, \frac{\overline{q}\, -\, T_\kappa \left( \overline{q} \right)}{\sqrt{\Delta_\kappa + \beta}}\, \Delta_\kappa \right|\, \leqslant\ f''(u^\ell)\, |\overline{q}|\, \sqrt{\Delta_\kappa}\ \leqslant\ \half\, f''(u^\ell) \left( \overline{q}^2\, +\, \Delta_\kappa \right)\, \leqslant\ \half\, \|f''(u^\ell)\|_{L^\infty} \left( \overline{q}^2\, +\, \half\, \overline{q^2} \right).
\end{equation*}
Using \eqref{TS} we obtain 
\begin{equation*}
\left| \half\, \widetilde{P}\, \frac{\overline{T_\kappa\left(q\right)}\, -\, T_\kappa\left(\overline{q}\right)}{\sqrt{\Delta_\kappa + \beta}} \right|\, 
\leqslant\ {\textstyle \frac{\sqrt{2}}{2}}\,\widetilde{P}.
\end{equation*}
Since the $L^1$ convergence implies the pointwise convergence (up to a subsequence), then, using the dominated convergence theorem with \eqref{uLinf}, \eqref{Fatou1} and the fact that $\widetilde{P} \in L^1_{loc}$, we obtain
\begin{equation*}
\lim_{\kappa \to \infty}
\left\| \half\, f''(u^\ell)\, \frac{\overline{q}\, -\, T_\kappa \left( \overline{q} \right)}{\sqrt{\Delta_\kappa + \beta}}\, \Delta_\kappa \right\|_{L^1(\Omega)}\, +\ \lim_{\kappa \to \infty} \left\| \half\, \widetilde{P}\, \frac{\overline{T_\kappa\left(q\right)}\, -\, T_\kappa\left(\overline{q}\right)}{\sqrt{\Delta_\kappa + \beta}} \right\|_{L^1(\Omega)}\, =\ 0.
\end{equation*}
for any compact set $\Omega \subset (0,\infty)\times \R$.
Taking $\kappa \to \infty$ in \eqref{Delta_k} we obtain 
\begin{equation}\nonumber
\left[\sqrt{\Delta + \beta} \right]_t\ +\ \left[ f'(u^\ell)\, \sqrt{\Delta + \beta} \right]_x\ \leqslant\ \frac{\beta\, \overline{q}\, f''(u^\ell)}{ \sqrt{\Delta + \beta}}\ \leqslant\ \sqrt{\beta}\, |\overline{q}|\, f''(u^\ell),
\end{equation}
where $\Delta \eqdef \overline{q^2}-\overline{q}^2$.
Taking now $\beta \to 0$ we obtain
\begin{equation}\label{Tg}
\left[\sqrt{\Delta} \right]_t\ +\ \left[ f'(u^\ell)\, \sqrt{\Delta } \right]_x\ \leqslant\ 0 \qquad \mathrm{in}\ (t_0,\infty) \times \R.
\end{equation}

\textbf{Step 5.} As in \cite{wave3}, let $\varphi \in C^\infty_c(\mathds{R})$ satisfying $\varphi(x)=1$ for $|x| \leqslant 1$ and $\varphi(x)=0$ for $|x| \geqslant 2$. Since $\sqrt{\Delta} \in L^\infty ((0,\infty), L^2(\mathds{R}))$, then, for all $n \geqslant 1$, we have $\sqrt{\Delta}  \varphi(x/n) \in L^\infty ((0,\infty), L^1(\mathds{R}))$. Then almost all $t>0$ are Lebesgue points of $t \mapsto \int_\mathds{R} \sqrt{\Delta}(t,x) \varphi(x/n) \mathrm{d}x,$ $\forall n \geqslant 1$.
Let $\bar{t}>0$ be a Lebesgue point of $t \mapsto \int_\mathds{R} \sqrt{\Delta} (t,x) \varphi(x/n) \mathrm{d}x$ and $\delta \in (0,\bar{t}/2)$. Let also $\psi \in C^\infty_c((0,\infty))$ satisfying
\begin{gather*}
\psi(t)\ =\ 0 \quad \mathrm{on}\quad (0,\delta/2) \cup (\bar{t}+\delta, \infty ), \qquad \psi(t)\ =\ 1 \quad \mathrm{on}\quad  (\delta, \bar{t}-\delta ),\\
0\ \leqslant\ \psi'(t)\ \leqslant\ C/\delta, \quad \mathrm{on}\quad (\delta/2,\delta), \qquad 
-\psi'(t)\ \geqslant\ C/\delta, \quad \mathrm{on}\quad (\bar{t}-\delta,\bar{t}+\delta).
\end{gather*}
Multiplying \eqref{Tg} by $\varphi(x/n) \psi(t)$, integrating on $(0,\infty) \times \mathds{R}$ and using integration by parts one obtains
\begin{gather*}
{\textstyle \frac{C}{\delta}}\, \int_{\bar{t}-\delta}^{\bar{t}+\delta} \int_\mathds{R} \sqrt{\Delta} (t,x)\, \varphi(x/n)\, \mathrm{d}x\, \mathrm{d}t\ 
\leqslant\ 
-\int_{\bar{t}-\delta}^{\bar{t}+\delta} \int_\mathds{R}  \sqrt{\Delta} (t,x)\, \varphi(x/n)\, \psi'(t)\, \mathrm{d}x\, \mathrm{d}t\\
\leqslant\ {\textstyle \frac{C}{\delta}}\, \int_{\delta/2}^{\delta} \int_\mathds{R} \sqrt{\Delta}(t,x)\, \varphi(x/n)\, \mathrm{d}x\, \mathrm{d}t\ +\
{\textstyle \frac{1}{n}}\, ||f'(u^\ell)||_{L^\infty}\, \int_{\delta/2}^{\bar{t}+\delta} \int_\mathds{R} \sqrt{\Delta}(t,x)\, \left| \varphi'(x/n) \right|\, \mathrm{d}x\, \mathrm{d}t.
\end{gather*}
From \eqref{tauto0}, we have
\begin{equation*}
\lim_{t \to 0} \int_\mathds{R} \sqrt{\Delta}(t,x)\, \varphi(x/n)\, \mathrm{d}x\ =\ 0 \quad \implies \quad 
\lim_{\delta \to 0} {\textstyle \frac{1}{\delta}}\, \int_{\delta/2}^{\delta}  \int_\mathds{R} \sqrt{\Delta}(t,x)\, \varphi(x/n)\, \mathrm{d}x\, \mathrm{d}t\ =\ 0.
\end{equation*}
Since $\bar{t}>0$ is a Lebesgue point of $t \mapsto \int_\mathds{R} \sqrt{\Delta}(t,x) \varphi(x/n) \mathrm{d}x$, then taking first $\delta \to 0$ and then $n \to \infty$ we obtain 
\begin{equation*}
\sqrt{\Delta}(\bar{t},x)\ =\ 0 \quad \mathrm{a.e.}\, (\bar{t},x) \in (0,\infty) \times \mathds{R}.
\end{equation*}
Hence $\overline{q^2} = \overline{q}^2$ almost everywhere, which implies that
$\mu_{t,x}^\ell(\xi) = \delta_{\overline{q}(t,x)}(\xi) = \delta_{u_x^\ell(t,x)}(\xi)$. 
\qed\\

\textit{Proof of Theorem \ref{thm:existence}}. All the limits in this proof are up to a subsequence.
Let $u^{\ell,\varepsilon}  $ be the solution given in Theorem \ref{Thm:existenceapp}. Then, from Lemma \ref{lemYoungm}, Lemma \ref{lem:2Dirac} and Lemma \ref{lem:alpha+2} we have that, as $\varepsilon \to 0$ 
\begin{align}\label{Wtconv_xi}
u_x^{\ell,\varepsilon} \qquad
&\rightharpoonup \qquad u_x^{\ell} \qquad
\mathrm{in}\ L^p_{loc}((0,\infty) \times \mathds{R}), \\ \nonumber
\left\|  \left(u_x^{\ell,\varepsilon}\right)^2 \right\|_{L^1(\Omega)}\quad
&\to \quad \left\|\left(u_x^{\ell}\right)^2\right\|_{L^1(\Omega)},
\end{align}
for any $p \in [2,3)$ and compact set $\Omega \subset (0,\infty) \times \mathds{R} $. This implies  that 
\begin{align}\label{Stconv_xi}
u_x^{\ell,\varepsilon} \qquad
&\to \qquad u_x^{\ell} 
&\mathrm{in}\ L^2_{loc}((0,\infty) \times \mathds{R}).
\end{align}
Using Lemma \ref{lem:alpha+2} and Lemma \ref{Strong_conv} we obtain that for all $p \in [2,3)$, we have 
\begin{align}\label{Wtconv_tau}
u_t^{\ell,\varepsilon}\  \qquad
&\rightharpoonup \qquad u_t^\ell
&\mathrm{in}\ L^p_{loc}((0,\infty) \times \mathds{R}).
\end{align}
From Lemma \ref{lem:2Dirac} and Lemma \ref{Strong_conv1.5} we obtain 
\begin{equation}\label{TildeP}
\widetilde{P}\ =\ \half\, \mathfrak{G}\ast \left\{ f''(u^\ell) \left(u^\ell_x\right)^{2}\right\}.
\end{equation}
Now, using Lemma \ref{Strong_conv1.5} and taking the weak limit $\varepsilon \to 0$ in \eqref{rBep} and \eqref{Ene_ep} we obtain \eqref{rB} and \eqref{Ene_local}.
Doing the proof of \eqref{tauto0} for any $t_0$, we obtain that \eqref{rcont}. 
The Oleinik inequality \eqref{Ol:main:thm} follows from Lemma \ref{Lem:Oleinik}.
The inequality \eqref{alpha+2_xi} follows from Lemma \ref{lem:alpha+2}, \eqref{Wtconv_xi} and \eqref{Wtconv_tau}.  
Finally, \eqref{TV} follows from Lemma \ref{lemTV} and  \eqref{Stconv_xi}. \qed

\begin{rem}
From \eqref{Sbar1}, Lemma \ref{lem:2Dirac} and \eqref{TildeP} we have 
\begin{gather}\nonumber
\left[S_\kappa\left(u_x^\ell \right)\right]_{t}\, +\, \left[f'(u^\ell)\, S_\kappa\left(u_x^\ell \right)\right]_{x}\ =\ -P\, T_\kappa\left(u_x^\ell\right)\\ \label{Sbar1.5}
+\ f''(u^\ell)\, u_x^\ell\, S_\kappa\left(u_x^\ell \right)\, -\ \half\, f''(u^\ell) \left(u_x^\ell\right)^2\, T_\kappa\left(u_x^\ell\right).
\end{gather}
This will be used in Section \ref{sec:infty} below.
\end{rem}

\section{The limiting case $\ell \to 0$}\label{sec:zero}

Let $u^\ell$ be the dissipative solutions of \eqref{rB} given by Theorem \ref{thm:existence}. The aim of this section is to study the limit of $u^\ell$ when $\ell \to 0$.

%

\subsection{Uniform estimates for small $\ell$}
Considering $\ell \leqslant 1$, then, from the energy equation \eqref{Ene_local} we have 
\begin{subequations}\label{lossenergy}
\begin{align}
\int_\mathds{R} \left( (u^\ell)^2\ +\ \ell^2\,  (u_x^\ell)^2\, \right) \mathrm{d}x\  &\leqslant\ \left\|u_0 \right\|_{H^1}^2,\\ \label{ell2P}
\int_\mathds{R} \ell^2\, P\, \mathrm{d}x\ =\ \int_\mathds{R} \half\, \ell^2\, f''(u^\ell)\, (u_x^\ell)^2\, \mathrm{d}x\  &\leqslant\ \half\, C\, \left\|u_0 \right\|_{H^1}^2.
\end{align}
\end{subequations}
In order to prove that $u^\ell$ converges to the unique entropy solution of the scalar conservation law \eqref{SCL} we need : $(i)$ a strong convergence of $u^\ell$ to some $u^0$; $(ii)$ the weak convergence of $\ell^2 P \rightharpoonup 0$. The result then follows directy by taking $\ell \to 0$ in \eqref{rB} and using \eqref{Ol:main:thm}. However, the estimate \eqref{ell2P} in not enough to identify the weak limit of $\ell^2 P$ as $\ell \to 0$. The following lemma gives a better estimation of $\ell^2 P$.

\begin{lem}\label{mu=0}
For all compact set $\Omega \subset (0,\infty) \times \mathds{R}$ and $\alpha \in (0,1)$ there exists $C_{\Omega,\alpha}>0$ such that for all $\ell \leqslant 1$ we have
\begin{align}\label{e}
\int_\Omega \ell^2\, P\, \mathrm{d}x\, \mathrm{d}t\ &\leqslant\ \ell^{\frac{2\alpha}{2+\alpha}}\, C_{\Omega,\alpha}.
\end{align}
\end{lem}
\proof
Without losing of generality, we take $\Omega = [t_1,t_2] \times [a,b]$ such that $t_1>0$. We take also $\alpha = 2k/(2k+1)$ such that $k \in \mathds{N}$.
Let $\psi \in C^\infty_c(]0,\infty[ \times \mathds{R})$ such that $\psi \geqslant 0$, $\psi=1$ on $\Omega$ and $\mathrm{supp}(\psi) \subset [t_1/2,t_2+1] \times [a-1,b+1]$.

\textbf{Step 1.} Multiplying \eqref{rBep_x} by $\ell^2 |q^{\ell,\varepsilon}|$ we obtain
\begin{equation*}
\left[\half\, \ell^2\, q^{\ell,\varepsilon}\, |q^{\ell,\varepsilon}| \right]_t\ +\ \left[\half\, \ell^2\, f'(u^{\ell,\varepsilon})\, q^{\ell,\varepsilon}\, |q^{\ell,\varepsilon}| \right]_x\ +\ \ell^2\, |q^{\ell,\varepsilon}|\, P^\varepsilon\ -\ \half\, \ell^2\, f''(u^{\ell,\varepsilon})\, |q^{\ell,\varepsilon}|\, \chi_\varepsilon(q^{\ell,\varepsilon})\ =\ 0.
\end{equation*}
Multiplying by $\psi$ and integrating we obtain
\begin{align*}
\int_\Omega \ell^2\, |q^{\ell,\varepsilon}|\, P^\varepsilon\, \mathrm{d}x\, \mathrm{d}t\ 
\leqslant &\ \int_{(0,\infty) \times \mathds{R}} \ell^2\, |q^{\ell,\varepsilon}|\, P^\varepsilon\, \psi\, \mathrm{d}x\, \mathrm{d}t\\
=&\ \half\, \ell^2\, \int_{(0,\infty) \times \mathds{R}} \left[ \psi_t\,  q^{\ell,\varepsilon}\, |q^{\ell,\varepsilon}|\ +\ \psi_x\, f'(u^{\ell,\varepsilon})\, q^{\ell,\varepsilon}\, |q^{\ell,\varepsilon}|\right] \mathrm{d}x\, \mathrm{d}t\\
&+\ \half\, \ell^2\, \int_{(0,\infty) \times \mathds{R}} \psi\, f''(u^{\ell,\varepsilon})\, |q^{\ell,\varepsilon}|\, \chi_\varepsilon(q^{\ell,\varepsilon})\, \mathrm{d}x\, \mathrm{d}t.
\end{align*}
Then, from \eqref{eneequation} and \eqref{TVep} we have
\begin{equation}\label{a}
\int_\Omega \ell^2\, |q^{\ell,\varepsilon}|\, P^\varepsilon\, \mathrm{d}x\, \mathrm{d}t\ \leqslant\ C_{\Omega,\alpha}, \quad \forall \varepsilon>0, \forall \ell \in  (0,1).
\end{equation}

\textbf{Step 2.} 
Using \eqref{a} and \eqref{lossenergy} we have 
\begin{align}\nonumber
\int_\Omega \ell^2\, |q^{\ell,\varepsilon}|^\alpha\, P^\varepsilon\, \mathrm{d}x\, \mathrm{d}t\
&=\ \int_{\Omega \cap \{ |q^{\ell,\varepsilon}| > 1 \} } \ell^2\, |q^{\ell,\varepsilon}|^\alpha\, P^\varepsilon\, \mathrm{d}x\, \mathrm{d}t\ +\ \int_{\Omega \cap \{ |q^{\ell,\varepsilon}| \leqslant 1 \} } \ell^2\, |q^{\ell,\varepsilon}|^\alpha\, P^\varepsilon\, \mathrm{d}x\, \mathrm{d}t\\ \nonumber
&\leqslant\ \int_{\Omega \cap \{ |q^{\ell,\varepsilon}| > 1 \} } \ell^2\, |q^{\ell,\varepsilon}|\, P^\varepsilon\, \mathrm{d}x\, \mathrm{d}t\ +\ \int_{\Omega \cap \{ |q^{\ell,\varepsilon}| \leqslant 1 \} } \ell^2\, P^\varepsilon\, \mathrm{d}x\, \mathrm{d}t\\ \label{b}
&\leqslant\ C_{\Omega,\alpha}, \quad \forall \varepsilon>0, \forall \ell \in  (0,1).
\end{align}

\textbf{Step 3.} Multiplying \eqref{rBep_x} by $\ell^2 (q^{\ell,\varepsilon})^\alpha$ we obtain
\begin{align*}
\frac{1-\alpha}{2\, ( \alpha + 1)}\, \ell^2\, f''(u^{\ell,\varepsilon}) (q^{\ell,\varepsilon})^{2+\alpha}\
=&\ \ell^2\, \left( \frac{(q^{\ell,\varepsilon})^{1+\alpha}}{1+\alpha} \right)_t\ +\ \ell^2\, \left( \frac{f'(u^{\ell,\varepsilon})\, (q^{\ell,\varepsilon})^{1+\alpha}}{1+\alpha} \right)_x\ +\ \ell^2\, (q^{\ell,\varepsilon})^\alpha\, P^\varepsilon\\
&-\ \half\, \ell^2\, f''(u^{\ell,\varepsilon})\, (q^{\ell,\varepsilon})^\alpha\, \chi_\varepsilon(q^{\ell,\varepsilon})\\
\leqslant &\ \ell^2\, \left( \frac{(q^{\ell,\varepsilon})^{1+\alpha}}{1+\alpha} \right)_t\ +\ \ell^2\, \left( \frac{f'(u^{\ell,\varepsilon})\, (q^{\ell,\varepsilon})^{1+\alpha}}{1+\alpha} \right)_x\ +\ \ell^2\, (q^{\ell,\varepsilon})^\alpha\, P^\varepsilon.
\end{align*}
Multiplying by $\psi$, doing as in step 1 and using \eqref{b} we obtain 
\begin{equation*}
\int_\Omega \ell^2\, (q^{\ell,\varepsilon})^{2+\alpha}\, \mathrm{d}x\, \mathrm{d}t\ \leqslant\ C_{\Omega,\alpha}, \quad \forall \varepsilon>0, \forall \ell \in  (0,1).
\end{equation*}
Using now Holder's inequality, we obtain 
\begin{equation}\label{d}
\int_\Omega \ell^2\, (q^{\ell,\varepsilon})^{2}\, \mathrm{d}x\, \mathrm{d}t\ \leqslant\ \ell^{\frac{2\alpha}{2+\alpha}}\, C_{\Omega,\alpha}, \quad \forall \varepsilon>0, \forall \ell \in  (0,1).
\end{equation}

\textbf{Step 4.} The equation \eqref{rBep_x} can be rewritten as
\begin{equation}\label{rBep_x2}
q^{\ell,\varepsilon}_{t}\, +\, \left[ f'(u^{\ell,\varepsilon})\, q^{\ell,\varepsilon}\right]_{x}\ +\ P^\varepsilon\ -\ \half\, f''(u^{\ell,\varepsilon})\, (q^{\ell,\varepsilon})^2\  -\ \half\, f''(u^{\ell,\varepsilon})\, \chi_\varepsilon(q^{\ell,\varepsilon})\ =\ 0.
\end{equation}
Multiplying \eqref{rBep_x2} by $\ell^2 \psi$ and doing as in step 1 and using \eqref{TVep} we obtain 
\begin{equation*}
\int_\Omega \ell^2\, P^\varepsilon\, \mathrm{d}x\, \mathrm{d}t\ \leqslant\ C_{\Omega,\alpha} \left[ \ell^2\, +\,  \int_\Omega \ell^2\, (q^{\ell,\varepsilon})^{2}\, \mathrm{d}x\, \mathrm{d}t \right]\ \leqslant\ \ell^{\frac{2\alpha}{2+\alpha}}\, C_{\Omega,\alpha}, \quad \forall \varepsilon>0, \forall \ell \in  (0,1].
\end{equation*}
Then \eqref{e} follows by taking $\varepsilon \to 0$ with Lemma \ref{Strong_conv1.5} and \eqref{TildeP}.
\qed

\subsection{Precompactness}

In order to obtain a strong compactness of $(u^\ell)_\ell$,  we use the Aubin--Lions--Simon theorem. For that purpose, we let $\mathcal{I} \subset \mathds{R}$ to be a compact interval and we define 
\begin{equation}\label{Wdef}
W(\mathcal{I})\ \eqdef\ \left\{ g\, \in\, \mathcal{D}'(\mathcal{I}),\ \exists\,
G\, \in\, L^1(\mathcal{I})\ \textrm{such that } G'\ =\ g \right\} ,
\end{equation}
where the norm of the space $W(\mathcal{I})$ is given by
\begin{equation}\label{norm}
\|g\|_{W(\mathcal{I})}\ \eqdef\ \inf_{c\, \in\, \mathds{R}} 
\|G + c\,\|_{L^1(\mathcal{I})}\ =\ \min_{c\, \in\, \mathds{R}} 
\|G + c\,\|_{L^1(\mathcal{I})} .
\end{equation}
\begin{lem} The space $W(\mathcal{I})$ is a Banach space and the 
embedding
\begin{equation}
L^1(\mathcal{I})\ \hookrightarrow\ W(\mathcal{I}),
\end{equation}
is continuous.
\end{lem}
\proof Let $(g_n)_{n\in\mathds{N}}$ be a Cauchy sequence in $W(\mathcal{I})$ and let 
$G_n$ be a primitive of $g_n$. From the definition of the norm \eqref{norm}, there exists 
a constant $c_n$ such that $(\tilde{G}_n-c_n)_{n\in\mathds{N}}$ (where $\tilde{G}_n=G_n
+c_n$) is a Cauchy sequence in $L^1(\mathcal{I})$. 
Let $\tilde{G}$ be the limit of $\tilde{G}_n$ in $L^1(\mathcal{I})$. Then
\begin{equation}
\|g_n\ -\ \tilde{G}'\|_{W(\mathcal{I})}\ \leqslant\ 
\|\tilde{G_n}\ -\ \tilde{G} \|_{L^1(\mathcal{I})},
\end{equation}
implying that $W(\mathcal{I})$ is a Banach space.

If $g \in L^1(\mathcal{I})$, then $G(x) - G(a) = \int_a^x g(y)\, \mathrm{d}y$  
for almost all $x,a \in \mathcal{I}$. Therefore, 
\begin{equation}
\|g\|_{W(\mathcal{I})}\ \leqslant\ \int_\mathcal{I} |G(x)\ -\ G(a)|\ 
\mathrm{d}x\ \leqslant\ |\mathcal{I}| \int_\mathcal{I} |g(y)|\ \mathrm{d}y 
,
\end{equation}
which ends the proof of the continuous embedding. \qed\\

The previous lemma and Helly's selection theorem imply that
\begin{equation}\label{embeddings}
W^{1,1}(\mathcal{I})\ \Subset\ L^1(\mathcal{I})\ \hookrightarrow\ 
W(\mathcal{I}),
\end{equation}
where the first embedding is compact and the second one is continuous.

The estimates \eqref{lossenergy} imply that $u^\ell$ is uniformly bounded on $L^\infty (\mathds{R}^+, 
L^2(\mathds{R}))$. Subsequently, it is also uniformly bounded on $L^\infty (\mathds{R}^+, 
L^1(\mathcal{I}))$. Then \eqref{TV} implies that $u^\ell$ is bounded on $L^\infty([0,T],
W^{1,1}(\mathcal{I}))$.

Since $0<c \leqslant f''(u) \leqslant C$, then $c u^2/2 \leqslant f(u)-f'(0) u -f(0) \leqslant C u^2/2$. Integrating \eqref{rB}, we obtain that $\int u^\ell\, \ud x = \int u_0\, \ud x $. This implies with  \eqref{lossenergy} that $f(u^\ell) - f(0)
+ \ell^2 P$ is uniformly bounded on $L^\infty([0,T],L^1(\mathcal{I}))$.  
Since 
$u^\ell_t = -\left( f(u^\ell) - f(0) + \ell^2 P \right)_x$, \eqref{norm} implies 
that $u^\ell_t$ is bounded on $L^\infty([0,T],W(\mathcal{I}))$. 
Then, using \eqref{embeddings} with Aubin--Lions--Simon lemma \cite{Simon}, we obtain that, up to a subsequence, $u^\ell$ converges to some $u^0$ in $C([0,T],L^1(\mathcal{I}))$ as $\ell \to 0$. Using an interpolation with \eqref{TV}, we obtain the convergence in $L^\infty([0,T],L^p(\mathcal{I}))$ for any $p\in [1,\infty)$.

The proof of Theorem \ref{thm:0} follows directly by taking $\ell \to 0$ in the weak formulation of \eqref{rB} and using Lemma \ref{mu=0} with the Oleinik inequality \eqref{Ol:main:thm}.
Due to the uniqueness of the entropy solution of the scalar conservation laws, we deduce that all the sequence $u^\ell$ converges to $u^0$.

\section{The limiting case $\ell \to \infty$}\label{sec:infty}

We consider $u^\ell$ the dissipative solutions of \eqref{rB} given by Theorem \ref{thm:existence}. The aim of this section is to study the limit of $u^\ell$ when $\ell \to \infty$.
 In Section \ref{sec:precompactness} above, a proof of the limiting case $\varepsilon \to 0$ is presented. In the present section we will use very similar techniques to establish the limiting case $\ell \to \infty$.
Recognizing the similarity between the two sections, we will streamline the explanations in this section, as they are closely parallel to the arguments outlined in Section \ref{sec:precompactness}.
For a complete understanding of the proof techniques in this section, we encourage the reader to refer back to Section \ref{sec:precompactness} when needed.

We start by obtaining some uniform estimates of $u^\ell$ when $\ell$ is far from $0$.
Considering $\ell \geqslant 1$, then the energy equation \eqref{Ene_local} implies 
\begin{equation}\label{lossenergy2}
	\int_\mathds{R} (u_x^\ell)^2\, \mathrm{d}x\  \leqslant\ \left\|u_0 \right\|_{H^1}^2.
\end{equation}
This implies that 
\begin{equation}\label{Linfbound}
\|u^\ell\|_{L^\infty_{loc}([0,\infty) \times \R)}\ +\ \|u^\ell\|_{L^\infty_{loc}([0,\infty),H^1_{loc}(\R)}\ \leqslant\ C.
\end{equation}
Using \eqref{Pdef} and Young inequality we obtain for all $p \in [1,\infty]$ that
\begin{equation}\label{PPx0.5} 
\|P\|_{L^p}\ \leqslant\ {\textstyle \frac{C}{2}}\, \|\mathfrak{G}\|_{L^p}\, \left\|u_x^\ell \right\|_{L^2}^2,
\quad \quad
\|P_x\|_{L^p}\ \leqslant\ {\textstyle \frac{C}{2}}\, \|\mathfrak{G}_x\|_{L^p}\, \left\|u_x^\ell \right\|_{L^2}^2.
\end{equation}
Using that $\ell \geqslant 1$ and \eqref{lossenergy2} we obtain 
\begin{equation}\label{PPx}
\ell\, \|P\|_{L^\infty}\ +\  \|P\|_{L^p}\ +\ \|P_x\|_{L^p}\ \leqslant\ C.
\end{equation}

\begin{lem}\label{lem:alpha+22}
Let $\alpha \in (0,1)$, $T>0$ and $ [a,b] \subset \R$, then there exists a constant $C=C(\alpha,T,a,b)>0$, such that for all $\ell \geqslant 1$ we have
\begin{equation}\label{alpha+22}
\int_0^T \int_a^b \left[ |u^{\ell}_t|^{2+\alpha}\ +\ |u^{\ell}_x|^{2+\alpha} \right] \mathrm{d}x\, \mathrm{d}t\ \leqslant\ C.
\end{equation}
\end{lem}
\proof When $\ell \geqslant 1$, one can use \eqref{lossenergy2}, \eqref{Linfbound} with \eqref{PPx0.5} and do the same proof of \eqref{alpha+2} to obtain a constant $C>0$ that does not depend on $\ell$. We conclude by taking $\varepsilon \to 0$. \qed

\begin{lem}\label{Strong_conv2} Let $[a,b] \subset \mathds{R}$ be a compact interval. Then, there exist $u^\infty \in L^\infty ([0,\infty), H^1([a,b]))$ and a subsequence of $(u^{\ell})_\ell$ such that, as $\ell \to \infty$, we have
\begin{align*}
u^{\ell} \quad &\to \quad u^\infty \qquad \mathrm{in}\  L^{\infty}([0,T] \times [a,b]),\ \forall T>0, \\
u^{\ell} \quad &\rightharpoonup \quad u^\infty \qquad \mathrm{in}\   H^1([0,T]\times [a,b]),\ \forall T>0.
\end{align*}
\end{lem}
\proof 
Using \eqref{Linfbound}, \eqref{rB} and \eqref{PPx} we obtain that 
\begin{equation}\label{L2L22}
\left\|u_t^{\ell}  \right\|_{L^2([0,T] \times [a,b])}\ \leqslant\ C_{T}.
\end{equation}
The weak convergence in $H^1([0,T] \times [a,b])$ follows directly. Using the inequality 
\begin{equation*}
\left\| u^{\ell}(t,\cdot)\, -\, u^{\ell}(s,\cdot) \right\|_{L^2([a,b])}^2\ =\ \int_a^b \left( \int_s^t u^{\ell}_t(\tau,x)\, \mathrm{d}\tau \right)^2 \mathrm{d}x\ \leqslant\ |t-s|\, \left\|u^{\ell}_t \right\|_{L^2([0,T] \times [a,b])}^2,
\end{equation*}
with \eqref{L2L22} we obtain that 
\begin{equation*}
\lim_{t \to s} \left\| u^{\ell}(t,\cdot)\, -\, u^{\ell} (s,\cdot) \right\|_{L^2([a,b])}\ =\ 0
\end{equation*}
uniformly on $\ell$. Then, using Theorem 5 in \cite{Simon} we can prove that that up to a subsequence,
$u^{\ell}$ converges uniformly to $u^\infty$ on any compact set of the form $[0,T] \times [a,b] $ as $\ell \to \infty$.
\qed

\begin{lem}\label{lemYoungm2}
There exist a subsequence of $\left\{q^{\ell} \right\}_\ell$ denoted also $\left\{q^{\ell} \right\}_\ell$ and a family of probability Young measures $\nu_{t,x}$ on $\mathds{R}$, such that for all functions $g \in C(\mathds{R})$ with $g(\xi) = \mathcal{O} (|\xi|^{2})$ at infinity, and for all $\varphi \in C^\infty_c((0,\infty) \times \mathds{R})$ we have
\begin{equation}\label{convegrgenceg2}
\lim_{\ell \to \infty}\, \int_{\mathds{R}^+ \times \mathds{R}} \varphi(t,x)\, g(q^{\ell})\, \mathrm{d}x\, \mathrm{d}t\ 
=\ \int_{\mathds{R}^+ \times \mathds{R}} \varphi(t,x) \int_{\mathds{R}} g(\xi)\,  \mathrm{d} \nu_{t,x}(\xi)\, \mathrm{d}x\, \mathrm{d}t.
\end{equation}
Moreover, the map
\begin{equation}\label{Fatou12}
(t,x)\ \mapsto\  \int_\mathds{R} \xi^2\, \mathrm{d} \nu_{t,x} (\xi)
\end{equation}
belongs to $L^\infty(\mathds{R}^+, L^1(\mathds{R})) $.
\end{lem}
\proof If $g(\xi)=\smallO(|\xi|^2)$, then the result is a direct consequence of \eqref{lossenergy2} and  Lemma \ref{lem:Young}.
If $g(\xi)=\mathcal{O}(|\xi|^2)$, let $\psi$ be a smooth cut-off function with $\psi(\xi)=1$ for $|\xi| \leqslant 1$ and $\psi(\xi)=0$ for $|\xi| \geqslant 2$, then 
\begin{equation}\label{convegrgencegpsi2}
\lim_{\ell \to \infty}\, \int_{\mathds{R}^+ \times \mathds{R}} \varphi(t,x)\, g_k(q^{\ell})\,  \mathrm{d}x\, \mathrm{d}t\ 
=\ \int_{\mathds{R}^+ \times \mathds{R}} \varphi(t,x) \int_{\mathds{R}} g_\kappa(\xi)\,  \mathrm{d} \nu_{t,x}(\xi,\zeta)\, \mathrm{d}x\, \mathrm{d}t,
\end{equation}
where $g_\kappa(\xi) \eqdef g(\xi) \psi\left( {\textstyle \frac{\xi}{\kappa}}\right)$ with $\kappa>0$. Using Holder inequality, Lemma \ref{lem:alpha+22} and $\Omega = \mathrm{supp}(\varphi)$ we obtain
\begin{gather*}
\left| \int_{\mathds{R}^+ \times \mathds{R}} \varphi(t,x)\, \left( g(q^{\ell})\, -\, g_\kappa(q^{\ell}) \right)\, \mathrm{d}x\, \mathrm{d}t \right|\
\leqslant\ 
\int_{\mathrm{supp}(\varphi) \cap \{|q^{\ell}| \geqslant \kappa\}} |\varphi(t,x)|\, \left|g(q^{\ell})\right|\, \mathrm{d}x\, \mathrm{d}t \\
\leqslant\ C \left( \int_{\mathrm{supp}(\varphi)} \left|g(q^{\ell})\right|^{p/2}\, \mathrm{d}x\, \mathrm{d}t \right)^{2/p} \left( \int_{\mathrm{supp}(\varphi) \cap \{|q^{\ell}| \geqslant \kappa\}} \mathrm{d}x\, \mathrm{d}t \right)^\frac{p-2}{p}
\\
\leqslant\ C\, \Big[ \big|\! \left\{(t,x) \in \mathrm{supp}(\varphi), |q^{\ell}| \geqslant \kappa \right\}\! \big| \Big]^\frac{p-2}{p}\
\leqslant\ C\, \kappa^{2-p}.
\end{gather*}
where $2<p<3$.
The last inequality with \eqref{convegrgencegpsi2} imply that we can interchange the limits $\kappa \to \infty$ and $\ell \to \infty$. Using that $|g_\kappa| \leqslant |g|$ and the dominated convergence theorem we obtain \eqref{convegrgenceg2}.
\qed

For the sake of simplicity, if no confusion with \eqref{bardef} is caused, we define in this section 
\begin{equation}\label{bardef2}
\overline{g(q)}\ \eqdef\ \int_{\mathds{R}} g(\xi)\, \mathrm{d} \nu_{t,x}(\xi)
\end{equation}
which is from \eqref{convegrgenceg2} the weak limit of $g(q^{\ell})$ as $\ell \to \infty$.

\begin{lem} As $t \to 0$ we have
\begin{gather} \label{tauto02}
\left\|u^\infty-u_0 \right\|_{\dot{H}^1(\mathds{R})}\ \to\ 0 \quad \qquad \mathrm{and} \qquad \quad
\int_\mathds{R} \left( \overline{q^2}\ -\ \overline{q}^2\, \right) \mathrm{d}x\ \to\ 0.
\end{gather}
\end{lem}

\proof

From Lemma \ref{Strong_conv2} and Lemma \ref{lemYoungm2} we have for all $t > 0$
\begin{gather*}
q^\ell(t,\cdot)\ \rightharpoonup\ \overline{q}(t,\cdot)\, =\, u^\infty_x(t,\cdot) \qquad \mathrm{as}\ \ell\, \to\, \infty \qquad \mathrm{in\ } L^2(\mathds{R}),\\
q^{\ell}(t,\cdot)^2\ \rightharpoonup\ \overline{q^2}(t,\cdot) \qquad \mathrm{when}\ \ell\, \to\, \infty \qquad \mathrm{in\ } \mathcal{D}'(\mathds{R}).
\end{gather*}
This, with Jensen's inequality and the energy equation \eqref{Ene_local} imply that 
\begin{align}\nonumber
\left\|u^\infty_x(t) \right\|_{L^2(\mathds{R})}^2\ &\leqslant\ \left\| \sqrt{\overline{q^2}}(t) \right\|_{L^2(\mathds{R})}^2\  \leqslant\
\liminf_{\ell \to \infty} \left\| u^{\ell}_x(t)\right\|_{L^2(\mathds{R})}^2\\ \label{Energyfin2}
& \leqslant\ \lim_{\ell \to \infty} \left\|\left(\ell^{-1}\, u_0,  u_0' \right)\right\|_{L^2(\mathds{R})}^2\ =\ \left\| u'_0 \right\|_{L^2(\mathds{R})}^2.
\end{align}
The inequality \eqref{lossenergy2} implies that the sequence $(q^{\ell})_\ell$ is bounded in the space $L^\infty([0,\infty), L^2(\mathds{R}))$.
Using \eqref{barqeqell}, we can prove that for all $T>0$ and for all $\varphi \in H^1(\mathds{R})$, the map 
\begin{equation}
t\ \mapsto\ \int_\mathds{R} \varphi(x)\,  q^{\ell}\, \mathrm{d}x
\end{equation}
is uniformly (on $t \in [0,T]$ and $\ell \geqslant 1 $) continuous. Then, Lemma \ref{lem:C_w} implies that 
\begin{equation}\label{Wconv22}
u_x^\infty(t,\cdot)\ \rightharpoonup\ u_0'\qquad \mathrm{when}\ t\, \to\, 0 \qquad \mathrm{in\ } L^2(\mathds{R}),
\end{equation}
which implies that 
\begin{equation*}
\left\|u_0' \right\|_{L^2}^2\ \leqslant\ \liminf_{t \to 0}  \left\| u^\infty_x(t)\right\|_{L^2}^2.
\end{equation*}
On another hand, \eqref{Energyfin2} implies
\begin{equation*}
\limsup_{t \to 0} \left\| u^\infty_x(t)\right\|_{L^2}^2\ \leqslant\ \left\|u_0' \right\|_{L^2}^2,
\end{equation*}
then 
\begin{equation}\label{Normconv3}
\lim_{t \to 0} \left\| u^\infty_x(t)\right\|_{L^2}^2\ =\ \left\|u_0' \right\|_{L^2}^2,
\end{equation}
which implies with \eqref{Wconv22} that 
\begin{equation}\label{Sconv2}
u_x^\infty(t,\cdot)\ \to\ u_0'\qquad \mathrm{when}\ t\, \to\, 0 \qquad \mathrm{in\ } L^2(\mathds{R}),
\end{equation}
The inequality \eqref{Energyfin2} with \eqref{Normconv3} imply 
\begin{equation}\label{Normconv22}
\lim_{t \to 0}\left\|\sqrt{\overline{q^2}} (t) \right\|_{L^2(\mathds{R})}^2\ =\
\lim_{t \to 0}\left\| u^\infty_x(t) \right\|_{L^2(\mathds{R})}^2\ =\ \left\| u'_0 \right\|_{L^2(\mathds{R})}^2.
\end{equation}
Then \eqref{tauto02} follows directly from \eqref{Sconv2} and \eqref{Normconv22}. 
\qed

\begin{lem}\label{lem:2Dirac2}
 The measure $\nu_{t,x}$ given in Lemma \ref{lemYoungm2} is a Dirac measure, and
\begin{equation}
\nu_{t,x}(\xi)\ =\ \delta_{u_x^\infty(t,x)}(\xi).
\end{equation}
\end{lem}
\proof 
\textbf{Step 0.} Let $S_\kappa$ and $T_\kappa$ defined as in \eqref{Sdef} and \eqref{Tdef} respectively.
As in Lemma \ref{lem:ST_k}, one can prove that for all $T>0$, we have
\begin{gather}\label{L1conv2}
\lim_{\kappa \to \infty} \left\| \overline{T_\kappa(q)}\ -\ T_\kappa\left(\overline{q}\right) \right\|_{L^1([0,T] \times \mathds{R})}\ =\
\lim_{\kappa \to \infty} \left\| \overline{T_\kappa(q)}\ -\ \overline{q} \right\|_{L^1([0,T] \times \mathds{R})}\ =\ 0. 
\end{gather}

\textbf{Step 1.}
Taking $\ell \to \infty$ in \eqref{Sbar1.5} and using \eqref{PPx}, we obtain 
\begin{equation}\label{Sq7.1}
\left[ \overline{S_\kappa \left(q\right)}\right]_{t}\, +\, \left[ f'(u^\infty)\, \overline{S_\kappa \left( q\right)} \right]_{x}\, =\ 
\half\, f''(u^\infty) \left[2\, \overline{q\, S_\kappa\left(q\right)}\, -\, \overline{q^2\, T_\kappa\left(q\right)}  \right].
\end{equation}
Taking $\ell \to \infty$ in \eqref{barqeqell}, we obtain 
\begin{gather}
\overline{q}_t\, +\, \left[f'(u^\infty)\,  \overline{q}\right]_{x}\  =\ \half\, f''(u^\infty)\, \overline{q^2}.
\end{gather}
Let $j_\varepsilon$ be a Friedrichs mollifier and $\overline{q}^\varepsilon \eqdef \overline{q} \ast j_\varepsilon$ then using Lemma \ref{lemRenor} we obtain
\begin{equation*}
\overline{q}^\varepsilon_t\, +\, \left[f'(u^\infty)\,  \overline{q}^\varepsilon\right]_{x}\  =\ \half\, f''(u^\infty)\, \overline{q^2}\ +\ \theta_\varepsilon.
\end{equation*}
where $\theta_\varepsilon \to 0$ in $L^1_{loc}([0,\infty) \times \R)$ as $\varepsilon \to 0$. 
Multiplying by $T_\kappa(\overline{q}^\varepsilon)$, we obtain 
\begin{align*}
\left[S_\kappa(\overline{q}^\varepsilon)\right]_t\, +\, \left[f'(u^\infty)\,  S_\kappa(\overline{q}^\varepsilon)\right]_{x}\,  =&\ - f''(u^\infty)\, \overline{q}\, \overline{q}^\varepsilon\, T_\kappa(\overline{q}^\varepsilon)\ +\   f''(u^\infty)\, \overline{q}\, S_\kappa(\overline{q}^\varepsilon)\\
&+\,  \left[ \half\, f''(u^\infty)\, \overline{q^2}\ +\ \theta_\varepsilon\right] T_\kappa(\overline{q}^\varepsilon).
\end{align*}
Taking $\varepsilon \to 0$, we obtain
\begin{align}\nonumber
\left[S_\kappa(\overline{q})\right]_t + \left[f'(u^\infty)\,  S_\kappa(\overline{q})\right]_{x}  =&\, - f''(u^\infty)\, \overline{q}^2\, T_\kappa(\overline{q})\,
+\, f''(u^\infty)\, \overline{q}\, S_\kappa(\overline{q})\,
+\, \half\, f''(u^\infty)\, \overline{q^2}\, T_\kappa(\overline{q})\\ \label{Sq7.2}
=&\ \half\, f''(u^\infty) \left[ T_\kappa(\overline{q}) \left(\overline{q^2}\, -\, \overline{q}^2 \right) +\, 2\, \overline{q}\, S_\kappa(\overline{q})\, -\,  \overline{q}^2\, T_\kappa(\overline{q}) \right].
\end{align}

\textbf{Step 2.} From \eqref{Sq7.1}, \eqref{Sq7.2}, \eqref{xiTS} and \eqref{Sdef} we obtain
\begin{gather}\nonumber
\left[\overline{S_\kappa \left(q\right)}\, -\, S_\kappa(\overline{q})\right]_t\, +\, \left[f'(u^\infty) \left( \overline{S_\kappa \left(q\right)}\, -\,  S_\kappa(\overline{q}) \right)\right]_{x}\, 
=\\  \nonumber
%
%
\half\, f''(u^\infty) \Big[
 \left( T_\kappa\left(\overline{q}\right)\, +\, \kappa \right)\, (\overline{q}\, +\, \kappa)^2\, \mathds{1}_{\overline{q} \leqslant -\kappa}\
+\ \left( T_\kappa\left(\overline{q}\right)\, -\, \kappa \right)\, (\overline{q}\, -\, \kappa)^2\, \mathds{1}_{\overline{q} \geqslant \kappa}  \\ \nonumber
 -\left( T_\kappa\left(\overline{q}\right)\, +\, \kappa \right)\, \overline{(q\, +\, \kappa)^2\, \mathds{1}_{q \leqslant -\kappa}}\
-\ \left( T_\kappa\left(\overline{q}\right)\, -\, \kappa \right)\, \overline{(q\, -\, \kappa)^2\, \mathds{1}_{q \geqslant \kappa}}\\ \label{Sq7.3}
-\ \kappa^2\, \left( \overline{T_\kappa(q)}\, -\, T_\kappa\left(\overline{q}\right) \right)\ 
-\ 2\, T_\kappa \left(\overline{q} \right) \left( \overline{S_\kappa(q)}\, -\, S_\kappa\left(\overline{q}\right) \right) \Big].
\end{gather}
From the definition \eqref{Tdef} we have 
\begin{equation}\label{7.1}
\left( T_\kappa\left(\overline{q}\right)\, +\, \kappa \right)\, (\overline{q}\, +\, \kappa)^2\, \mathds{1}_{\overline{q} \leqslant -\kappa}\
=\ \left( T_\kappa\left(\overline{q}\right)\, -\, \kappa \right)\, (\overline{q}\, -\, \kappa)^2\, \mathds{1}_{\overline{q} \geqslant \kappa}\ =\ 0.
\end{equation}
Since $T_\kappa\left(\overline{q}\right) \geqslant - \kappa $, then 
\begin{equation}\label{7.2}
 -\left( T_\kappa\left(\overline{q}\right)\, +\, \kappa \right)\, \overline{(q\, +\, \kappa)^2\, \mathds{1}_{q \leqslant -\kappa}}\ \leqslant\ 0.
\end{equation}
Let $t_0>0$ and $\kappa \geqslant 2/(c t_0)$, then from \eqref{Ol:main:thm}, we have for all $t \geqslant t_0$ that $q^\ell \leqslant \kappa $ and $\overline{q} \leqslant \kappa$. Then, using the convexity of $T_\kappa$ on $(-\infty,\kappa)$ and the Jensen's inequality we obtain 
\begin{equation}\label{7.3}
-\kappa^2\, \left( \overline{T_\kappa(q)}\, -\, T_\kappa\left(\overline{q}\right) \right)\ \leqslant\ 0, \qquad \forall t \geqslant t_0, \quad   \kappa \geqslant 2/(c\, t_0).
\end{equation}
We take again $t_0>0$ and $\kappa \geqslant 2/(c t_0)$, then for all $\varphi \in C^\infty_c((t_0,\infty) \times \mathds{R})$ we have 
\begin{equation}\label{7.4}
\int \overline{(q\, -\, \kappa)^2\, \mathds{1}_{q \geqslant \kappa}}\, \varphi\, \mathrm{d}x\, \mathrm{d}t\ =\ \lim_{\varepsilon \to 0} \int (q^\ell\, -\, \kappa)^2\, \mathds{1}_{q^\ell \geqslant \kappa}\, \varphi\, \mathrm{d}x\, \mathrm{d}t\ =\ 0.
\end{equation}
Defining $\tilde{\Delta}_\kappa \eqdef \overline{S_\kappa \left(q\right)}\, -\, S_\kappa\left(\overline{q}\right) $ and summing up \eqref{Sq7.3}, \eqref{7.1}, \eqref{7.2}, \eqref{7.3} and \eqref{7.4} we obtain 
\begin{equation*}
\left[\tilde{\Delta}_\kappa\right]_t\ +\ \left[ f'(u^\infty)\, \tilde{\Delta}_\kappa \right]_x\ \leqslant\ -\ f''(u^\infty)\, T_\kappa \left(\overline{q} \right) \tilde{\Delta}_\kappa .
\end{equation*}

\textbf{Step 3.}
Defining $\tilde{\Delta}^\varepsilon_\kappa \eqdef \tilde{\Delta}_\kappa \ast j_\varepsilon$ we obtain 
\begin{equation*}
\left[\tilde{\Delta}^\varepsilon_\kappa\right]_t\ +\ \left[ f'(u^\infty)\, \tilde{\Delta}^\varepsilon_\kappa \right]_x\ \leqslant\  -\ f''(u^\infty)\, T_\kappa \left(\overline{q} \right) \tilde{\Delta}_\kappa^\varepsilon\ +\ \tilde{\theta}_\varepsilon,
\end{equation*}
where $\tilde{\theta}_\varepsilon \to 0$ as $\varepsilon \to 0$ in $L^1_{loc}((0,\infty)\times \R)$.
Let $\beta >0$, multiplying by $\left( \tilde{\Delta}^\varepsilon_\kappa + \beta \right)^{-1/2}/2 $ we obtain 
\begin{gather*}\nonumber
\left[\sqrt{\tilde{\Delta}^\varepsilon_\kappa + \beta} \right]_t\ +\ \left[ f'(u^\infty)\, \sqrt{\tilde{\Delta}^\varepsilon_\kappa + \beta} \right]_x\ \leqslant\ \half\, f''(u^\infty)\, \frac{\overline{q}\, -\, T_\kappa \left( \overline{q} \right)}{\sqrt{\tilde{\Delta}^\varepsilon_\kappa + \beta}}\, \tilde{\Delta}^\varepsilon_\kappa\ +\ \frac{2\, \beta\, \overline{q}\, f''(u^\infty)\, +\, \tilde{\theta}_\varepsilon}{2\, \sqrt{\tilde{\Delta}^\varepsilon_\kappa + \beta}}.
\end{gather*}
Taking $\varepsilon \to 0$ we obtain 
\begin{gather}\label{Delta_k2}
\left[\sqrt{\tilde{\Delta}_\kappa + \beta} \right]_t\ +\ \left[ f'(u^\infty)\, \sqrt{\tilde{\Delta}_\kappa + \beta} \right]_x\ \leqslant\ \half\, f''(u^\infty)\, \frac{\overline{q}\, -\, T_\kappa \left( \overline{q} \right)}{\sqrt{\tilde{\Delta}_\kappa + \beta}}\, \tilde{\Delta}_\kappa\ +\ \frac{\beta\, \overline{q}\, f''(u^\infty)}{ \sqrt{\tilde{\Delta}_\kappa + \beta}}.
\end{gather}
Using that $|T_\kappa(\xi)| \leqslant |\xi|$ and $|S_\kappa(\xi)| \leqslant \xi^2/2$ we obtain 
\begin{align*}
\left| \half\, f''(u^\infty)\, \frac{\overline{q}\, -\, T_\kappa \left( \overline{q} \right)}{\sqrt{\tilde{\Delta}_\kappa + \beta}}\, \tilde{\Delta}_\kappa \right|\, &\leqslant\ f''(u^\infty)\, |\overline{q}|\, \sqrt{\tilde{\Delta}_\kappa}\ \leqslant\ \half\, f''(u^\infty) \left( \overline{q}^2\, +\, \tilde{\Delta}_\kappa \right)\\
 &\leqslant\ \half\, \|f''(u^\infty)\|_{L^\infty} \left( \overline{q}^2\, +\, \half\, \overline{q^2} \right)
\end{align*}
Since the $L^1$ convergence implies the pointwise convergence (up to a subsequence), then, using the dominated convergence theorem with Jensen inequality, \eqref{Linfbound} and \eqref{Fatou12}, we obtain
\begin{equation*}
\lim_{\kappa \to \infty}
\left\| \half\, f''(u^\infty)\, \frac{\overline{q}\, -\, T_\kappa \left( \overline{q} \right)}{\sqrt{\tilde{\Delta}_\kappa + \beta}}\, \tilde{\Delta}_\kappa \right\|_{L^1(\Omega)}\, =\ 0.
\end{equation*}
for any compact set $\Omega \subset (0,\infty)\times \R$.
Taking $\kappa \to \infty$ in \eqref{Delta_k2} we obtain 
\begin{equation}\nonumber
\left[\sqrt{\tilde{\Delta} + \beta} \right]_t\ +\ \left[ f'(u^\infty)\, \sqrt{\tilde{\Delta} + \beta} \right]_x\ \leqslant\ \frac{\beta\, \overline{q}\, f''(u^\infty)}{ \sqrt{\tilde{\Delta} + \beta}}\ \leqslant\ \sqrt{\beta}\, |\overline{q}|\, f''(u^\infty),
\end{equation}
where $\tilde{\Delta} \eqdef \overline{q^2}-\overline{q}^2$.
Taking now $\beta \to 0$ we obtain
\begin{equation}
\left[\sqrt{\tilde{\Delta}} \right]_t\ +\ \left[ f'(u^\infty)\, \sqrt{\tilde{\Delta} } \right]_x\ \leqslant\ 0 \qquad \mathrm{in}\ (t_0,\infty) \times \R.
\end{equation}

\textbf{Step 4.} Finaly, following Step 5 in the proof of Lemma \ref{lem:2Dirac} and using \eqref{tauto02}, we obtain that $\tilde{\Delta} = 0$ a.e.
\qed\\

\textit{Proof of Theorem \ref{thm:inf}}. All the limits in this proof are up to a subsequence.
Let $u^{\ell}  $ be the solution given in Theorem \ref{thm:existence}. Then, from Lemma \ref{lemYoungm2}, Lemma \ref{lem:2Dirac2}, Lemma \ref{lem:alpha+22} we have that, as $\ell \to \infty$
\begin{align}\label{Wtconv_xi2}
u_x^{\ell} \qquad
&\rightharpoonup \qquad u_x^{\infty} \qquad
\mathrm{in}\ L^p_{loc}([0,\infty) \times \mathds{R}), \\ \nonumber
\left\|  \left(u_x^{\ell}\right)^2 \right\|_{L^1(\Omega)}\quad
&\to \quad \left\|\left(u_x^{\infty}\right)^2\right\|_{L^1(\Omega)},
\end{align}
for any $p \in [2,3)$ and compact set $\Omega \subset [0,\infty) \times \mathds{R} $. This implies  that 
\begin{align}\label{Stconv_xi2}
u_x^{\ell} \qquad
&\to \qquad u_x^{\infty} 
&\mathrm{in}\ L^2_{loc}([0,\infty) \times \mathds{R}).
\end{align}
Using Lemma \ref{lem:alpha+22} and Lemma \ref{Strong_conv2} we obtain that for all $p \in [2,3)$, we have 
\begin{align}\label{Wtconv_tau2}
u_t^{\ell}\  \qquad
&\rightharpoonup \qquad u_t^\infty
&\mathrm{in}\ L^p_{loc}([0,\infty) \times \mathds{R}).
\end{align}
Now, taking the weak limit $\ell \to \infty$ in \eqref{barqeqell} we deduce that $u^\infty$ satisfies \eqref{gHS}. We multiply \eqref{Ene_local} by $\ell^{-2}$ and we take the limit $\ell \to \infty$ using \eqref{Linfbound} with \eqref{PPx} we obtain \eqref{Ene_localHS}.
Doing the proof of \eqref{tauto02} for any $t_0$, we obtain that \eqref{rcont2}. 
The Oleinik inequality \eqref{Ol:main:thm} implies $u_x^{\infty}(t,x)\ \leqslant\ \frac{1}{c\, t/2\, +\, 1/M}$ $a.e.\ (t,x) \in (0,\infty) \times \R$.
Using \eqref{alpha+22}, \eqref{Wtconv_xi2} and \eqref{Wtconv_tau2} we obtain that  $u^\infty_t,u^\infty_x \in L^{2+\alpha}_{loc}([0,\infty)\times \R)$, $\forall \alpha \in [0,1)$.
Finally, \eqref{TVgHS} follows from \eqref{TV} and \eqref{Stconv_xi2}. \qed

\appendix
\section{Some classical lemmas}\label{App:A}

Here, we recall simple versions of some classical lemmas that are needed in this paper. 

We start this section by the following lemma on the Young measures.
\begin{lem}(\cite{Focusing})\label{lem:Young}
Let $\mathscr{O}$ be a subset of $\mathds{R}^n$ with a zero-measure boundary. For any bounded family $\{v^\varepsilon \}_\varepsilon \subset L^p(\mathscr{O},\mathds{R}^N)$ with $p>1$ there exists a subsequence denoted also $\{v^\varepsilon \}_\varepsilon$ and a family of probability measures on $\mathds{R}^N$, $\left\{ \mu_{y}, y \in \mathscr{O} \right\}$ such that for all $f \in C^0(\mathds{R}^N)$ with $f(\xi)=\smallO(|\xi|^p)$ at infinity and for all $\phi \in C^\infty_c(\mathscr{O})$ we have 
\begin{equation}
\lim_{\varepsilon \to 0}\, \int_\mathscr{O} \phi(y)\, f(v^\varepsilon(y))\, \mathrm{d}y\ =\ \int_\mathscr{O} \phi(y) \int_\mathds{R} f(\xi)\, \mathrm{d} \mu_{y}(\xi)\, \mathrm{d}y
\end{equation}
with 
\begin{equation}\label{Fatou0}
\int_\mathscr{O} \int_\mathds{R} |\xi|^p\, \mathrm{d} \mu_y(\xi)\, \mathrm{d}y\ \leqslant\ \liminf_{\varepsilon \to 0} \|u^\varepsilon\|_{L^p(\mathscr{O})}^p.
\end{equation}
\end{lem}
Also, some other results on strong and weak precompactness are needed, then we recall.
\begin{lem}(\cite{Evans})\label{lem:fg}
Let $\Omega$ be an open set of $\mathds{R}^n$, assuming that $f_n \to f$ in $L^p(\Omega)$ with $p \in (1,\infty)$, $g_n$ is bounded in $L^q$ with $q \in (1,\infty)$ and $g_n \rightharpoonup  g $ in $L^q(\Omega)$, then for any $\varphi \in L^r(\Omega)$ such that $1/p+1/q+1/r=1$, we have
\begin{equation}
\lim_{n \to  \infty}\int_\Omega f_n\, g_n\, \varphi\, \mathrm{d}x\ =\ \int_\Omega f\, g\, \varphi\, \mathrm{d}x.
\end{equation}
\end{lem}

\begin{lem}(Lemma II.1 in \cite{DiPernaL1989})\label{lemRenor}
Let $c \in L^1_{loc}(\mathds{R}^+,H^1_{loc}(\mathds{R})$ and $f \in L^\infty_{loc}(\mathds{R}^+,L^2_{loc}(\mathds{R}))$. Let also $j_\varepsilon$ be a Friedrichs mollifier, then 
\begin{equation}
(c\, \partial_x f)\, \ast j_\varepsilon\ -\ c\, (\partial_x f \ast j_\varepsilon)\\\xrightarrow{\varepsilon \to 0} 0, \qquad \mathrm{in} \qquad L^1_{loc}(\mathds{R}^+ \times \mathds{R}).
\end{equation}
\end{lem}
\begin{lem}(Lemma C.1 in \cite{Lions})\label{lem:C_w}
Let $(f_n)_n$ be a bounded sequence in $L^\infty([0,T], L^2(\mathds{R}))$. If $f_n$ belongs to $C([0,T], H^{-1}(\mathds{R}))$ and  for any $\varphi \in H^1(\mathds{R})$, the map 
\begin{equation*}
t\ \mapsto\ \int_\mathds{R} \varphi(x)\, f_n(t,x)\, \mathrm{d}x
\end{equation*}
is uniformly continuous for $t \in [0,T]$ and $n \geqslant 1$, then $(f_n)_n$ is relatively compact in the space $C([0,T],L^2_w(\mathds{R}))$, where $L^2_w$ is the $L^2$ space equipped with its weak topology.
\end{lem}

\end{document}